\documentclass[11pt,a4paper]{article}
\usepackage[english,french]{babel}
\usepackage[latin1]{inputenc}
\usepackage[centertags]{amsmath}
\usepackage{amsfonts}
\usepackage{amssymb}
\usepackage{amsthm}
\usepackage{newlfont}
\usepackage{color}      
\usepackage{epsfig}

\input{psfig}

\usepackage{graphicx}

\usepackage{float}

\newtheorem{theo}{Theorem}
\newtheorem{prop}{Proposition}

\newtheorem{lem}{Lemma}

\newcommand{\complex}{\mathbb C}
\newcommand{\integer}{\mathbb N}
\newcommand{\rinteger}{\mathbb Z}
\newcommand{\real}{\mathbb R}

\newcommand{\eq}{\mathbb{E}_q}
\newcommand{\pochamer}[3]{\left(#1;#2\right)_{#3}}

\newcommand{\qhypermatricetcour}[3]{A(#1;#2;#3)}
\newcommand{\qhypertcour}[3]{\vphantom{}_3 \phi_2 \left(#1 ;  #2 ; #3 \right)}

\newcommand{\pmatrice}[4]{\begin{pmatrix}#1&#2\\#3&#4 \end{pmatrix}}

\newcommand{\vectt}[3]{\left(\begin{array}{*{100}c} #1 \\ #2 \\#3 \end{array}\right)}
\newcommand{\pmatricet}[9]{\begin{pmatrix}#1&#2&#3\\#4&#5&#6\\#7&#8&#9\end{pmatrix}}

\newcommand{\qcar}[1]{e_{q,#1}}

\newcommand{\ftz}[3]{F^{(0)}(#1;#2;#3)}
\newcommand{\ftinf}[3]{F^{(\infty)}(#1;#2;#3)}

\newcommand{\ytz}[3]{Y^{(0)}(#1;#2;#3)}
\newcommand{\ytinf}[3]{Y^{(\infty)}(#1;#2;#3)}

\newcommand{\jtz}[1]{J^{(0)}(#1)}
\newcommand{\jtinf}[1]{J^{(\infty)}(#1)}

\newcommand{\bon}{\Omega}
\newcommand{\und}[1]{\underline{#1}}
\newcommand{\sigq}{\sigma_q}

\newcommand{\psl}{\text{\textsf{PSl}}_2(\complex)}
\newcommand{\pgl}{\text{\textsf{PGl}}_2(\complex)}
\newcommand{\gl}[1]{\text{\textsf{Gl}}_#1(\complex)}
\renewcommand{\sl}[1]{\text{\textsf{Sl}}_#1(\complex)}
\newcommand{\G}{\text{\textsf{G}}}
\newcommand{\gz}{\text{\textsf{G}}^{\text{\textsf{0}}}}
\newcommand{\gzder}{\text{\textsf{G}}^{\text{\textsf{0,der}}}}
\newcommand{\hyp}[2]{\mathcal{H}_q(#1;#2)}
\newcommand{\perm}{\text{\textsf{Perm}}}
\newcommand{\permc}{\text{\textsf{Perm}}_{\complex^*}}
\newcommand{\gln}{\text{\textsf{Gl}}_n}

\begin{document}
\selectlanguage{english}
\begin{center}
\begin{huge}Galois groups of the Lie-irreducible generalized
$q$-hypergeometric equations of order three with $q$-real parameters :
an approach using a density theorem\end{huge}\\
$_{}$\\
\begin{LARGE} Julien Roques
 \end{LARGE}\\
$_{}$\\
 15th of November 2007\\
\end{center}
\hrule $_{}$ \vskip 20 pt

\noindent \textbf{Abstract}. \textit{In this paper we compute the
difference Galois groups of the Lie-irreducible regular singular
generalized $q$-hypergeometric equations of order 3 with $q$-real
parameters by using a density theorem due to Sauloy\footnote{An
alternative approach would be to use duality. This is a work in
progress and will appear elsewhere.}. In contrast with the
differential case, we show that these groups automatically contain
the special linear group $\sl{3}$.} \vskip 15 pt \noindent
\hrulefill \small \tableofcontents \normalsize
\newpage
In the whole paper $q$ is a complex number such that $0<|q|<1$.

\section{Generalized hypergeometric series and equations and $q$-ana\-logues}

\subsection{The differential case}

In accordance with the tradition, we denote by $\delta$ the Euler
differential operator : $\delta = z \frac{d~}{dz}$. Consider
$(r,s)\in\integer^{*2}$, $\und \alpha=(\alpha_1,...,\alpha_r)\in
\complex^r$ and $\und \beta=(\beta_1,...,\beta_s)\in \complex^s$.
Recall that the \textit{generalized hypergeometric series} with
parameters $(\und \alpha,\und \beta)$ is given by :
\begin{eqnarray*}
&&\sum_{n=0}^{+\infty} \frac{(\und \alpha)_n }{(\und \beta)_n}
z^n\\
&=&\sum_{n=0}^{+\infty} \frac{(\alpha_1)_n \cdots
(\alpha_r)_n}{(\beta_1)_n\cdots (\beta_s)_n} z^n\\
&=&\sum_{n=0}^{+\infty} \frac{\alpha_1 (\alpha_1+1) \cdots
(\alpha_1+n-1) \cdots \alpha_r (\alpha_r+1) \cdots (\alpha_r+n-1)}
{\beta_1 (\beta_1+1) \cdots (\beta_1+n-1) \cdots \beta_s
(\beta_s+1) \cdots (\beta_s+n-1)} z^n
\end{eqnarray*}
and that it provides us with a solution of the following linear
differential equation :
\begin{equation}\label{hyp class}\left(\prod_{j=1}^s(\delta + \beta_j - 1) -
z\prod_{i=1}^r(\delta + \alpha_i)\right)y=0\end{equation} called
the \textit{generalized hypergeometric equation} with parameters
$(\und \alpha,\und \beta)$.

These series and equations have retained the attention of many
authors. Since the pioneering work of Gauss, the theory of
generalized hypergeometric series and equations has figured
importantly in many branches of mathematics, from mathematical
physics to arithmetic. We shall now discuss more specifically the
Galois theory of hypergeometric equations.

As regards functional equations, the study of (\ref{hyp class})
leads us to distinguish two cases :
\begin{itemize}
\item[-] when $r=s$ the equation (\ref{hyp class}) is Fuchsian over the
Riemann  sphere with three singular points : $0$, $1$ and $\infty$;
\item[-] whereas when $r\neq s$ the equation (\ref{hyp class}) has two singular points~: $0$ and
$\infty$, one of these points being an \textit{irregular
singularity} (the other one is regular).
\end{itemize}
The Galois theory of the generalized hypergeometric equations in
both regular and irregular cases was in particular studied in
details by Beukers, Brownawell and Heckman in
\cite{beukersheckman,beukersbrownheck} and by Katz in
\cite{katzexpsum}. In \cite{duvalmitschi,mitschihyperconf} A.
Duval and C. Mitschi computed the Galois groups of some irregular
generalized hypergeometric equations by means of Ramis density
theorem.

From a technical point of view, in the Fuchsian case ($r=s$), the
determination of the Galois groups of the \textit{Lie-irreducible}
generalized hypergeometric equations (which are, roughly speaking,
equations with big Galois groups) rely on the fact that the local
monodromy of (\ref{hyp class}) around the singular point $z=1$ is
a pseudo-reflection : this allows us to apply a general result on
algebraic groups (see \cite{beukersheckman,katzexpsum}).

The present paper is concerned with the Galois theory of the
natural counterpart for $q$-difference equations of the
generalized hypergeometric equations with $r=s=3$.

\subsection{The $q$-difference case}

The natural $q$-analogues of the generalized hypergeometric
equations and series are respectively called the generalized
$q$-hypergeometric equations and series. Recall that the
\textit{generalized $q$-hypergeometric series} with parameters
$(\und{a},\und{b}) \in (\complex^*) ^r\times (\complex^*)^s$ is
defined by :
$$\sum_{n=0}^{+\infty} \frac{\pochamer{\und{a}}{q}{n} }{\pochamer{\und{b}}{q}{n} } z^n$$
and that, similarly to the differential case, it satisfies a
functional equation, namely the \textit{generalized
$q$-hypergeometric equation} with parameters $(\und{a},\und{b})$,
denoted by $\hyp{\und a}{\und b}$, and given by :
\begin{equation}\label{equa hypergeo} \left(\prod_{j=1}^s(\frac{b_j}{q} \sigq-1) -
z\prod_{i=1}^r(a_i \sigq-1)\right)\phi(z)=0\end{equation} where
$\sigq$ denotes the scaling operator acting on a function $\phi$
by $(\sigq\phi)(z)=\phi(qz)$.
We have used the classical notations for $q$-Pochhammer symbols :
$$\pochamer{a_i}{q}{n}=(1-a_i)(1-a_iq) \cdots (1-a_iq^{n-1})$$
and :
$$\pochamer{\und a}{q}{n}=\pochamer{a_1}{q}{n} \cdots
\pochamer{a_r}{q}{n}$$ (similar notations for $\und b$).
 It is usual to normalize both these equations and these series
by requiring $b_1=q$ : we will implicitly do this hypothesis in
the whole paper. The above $q$-hypergeometric series is then
denoted by $\vphantom{}_{r}\phi_{s-1}(\und a;\und b;z)$.

As in the differential case, one can easily check that the
$q$-difference equation $\hyp{\und a
}{\und b}$ is Fuchsian (see section \ref{section the basic objects} for this notion) if and only if $r=s$. \\

In this paper, we focus our attention on the case $r=s=3$.
The functional equation (\ref{hyp class}) is then equivalent to a
$3\times 3$ functional system. Indeed, with the following
notations :
$$\lambda(\und{a};\und{b};z)= \frac{1-z}{b_2b_3/q^2 - z a_1a_2a_3}, \ \ \ \
\mu(\und{a};\und{b};z)=\frac{z(a_1+a_2+a_3)-(1+b_2/q+b_3/q)}{b_2b_3/q^2
- z a_1a_2a_3}$$ and :
$$\delta(\und{a};\und{b};z)=\frac{(b_2b_3/q^2+b_2/q+b_3/q)-z(a_1a_2+a_2a_3+a_1a_3)}{b_2b_3/q^2
- z a_1a_2a_3}$$ a function $\phi$ is solution of the equation
$(\ref{equa hypergeo})$ if and only if the vector
$\Phi(z)=\vectt{\phi(z)}{\phi(qz)}{\phi(q^2z)}$ satisfies the
following functional system : \begin{equation}\label{syst
hypergeo} \Phi(qz) = \qhypermatricetcour{\und{a}}{\und{b}}{z}
\Phi(z) \end{equation} with :
$$\qhypermatricetcour{\und{a}}{\und{b}}{z}=\left( \begin{array}{ccc} 0&1&0\\0&0&1\\ \lambda(\und{a};\und{b};z) & \mu(\und{a};\und{b};z)&\delta(\und{a};\und{b};z) \end{array}\right).$$

Our main purpose is to show that the Galois groups of any
Lie-irreducible regular singular generalized $q$-hypergeometric
equation of order 3 with $q$-real parameters (that is of the form
$q^\alpha$ with $\alpha\in\real$) necessarily contains the special
linear group.
We emphasize that the similar statement for differential hypergeometric equations is \textit{false} (see \cite{katzexpsum}).\\


A number of authors have developed $q$-difference Galois theories
over the past years, among whom Franke \cite{frankepvdifference},
Etingof \cite{etingofgalois}, Van der Put and Singer
\cite{psgaloistheory}, Van der Put and Reversat \cite{prev},
Chatzidakis and Hrushovski \cite{chatzihru}, Sauloy
\cite{sauloyqgaloisfuchs}, Andr\'e \cite{andrenoncomm}, etc. The
exact relations between the existing Galois theories for
$q$-difference equations are partially understood. For this
question, we refer the reader to \cite{chatziharsing}, and also to
our Remark section \ref{constr}.

We will deal in this paper with the (Tannakian) Galois groups in
the sense of Sauloy in \cite{sauloyqgaloisfuchs} (see section
\ref{rappels galois} for an overview of the theory).\\ 

The notion of \textit{Lie-irreducibility} used
in this paper is a natural extension to $q$-difference equations
of a notion introduced by Katz : an equation is Lie-irreducible if
the linear representation of the neutral component of its Galois
group induced by the linear representation of the Galois group
itself given by the Tannakian duality is irreducible.\\

From a technical point of view, the fundamental fact which deeply
distinguishes the $q$-difference case from the differential case
is that we do not have presently a notion of local monodromy
around the singular point $z=1$ (or $q^{\rinteger}$);
nevertheless, we have some avatars built up from Birkhoff
connection matrices (generators of the connection component : see
section \ref{constr}). Whatever, it seems that we cannot use the
\textit{pseudo-reflection trick} mentioned above in order to
determine the Galois groups of the Lie-irreducible equations. This
is an important problem. Our method, in the three dimensional
case, is based on :
\begin{itemize}
\item[-] the classification of unimodular semi-simple algebraic subgroups of $\sl{3}$ acting irreducibly on $\complex^3$;
\item[-] a $q$-analogue of Schlesinger density theorem stated and established by Sauloy in \cite{sauloyqgaloisfuchs}.
 \end{itemize}

Last, we would like to point out Andr\'e's paper \cite{andrenoncomm}
which contains a nice computation of the Galois groups of some
generalized $q$-hypergeometric equations using a
\textit{specialization} procedure.\\


\subsection{Organization - Main result}\label{org}

The paper is organized as follows. In section \ref{rappels
galois}, we give a brief overview of the Galois theory of Fuchsian
$q$-difference systems in the sense of Sauloy. In section \ref{lie
irr} we introduce the notions of \textit{irreducibility}, of
\textit{Lie-irreducibility} and we discuss the irreducibility of the
generalized $q$-hypergeometric equations. In section \ref{pre alg
gps}, we investigate the properties of some algebraic subgroups of
$\gl{3}$. Section \ref{calcul} is devoted to the proof of our main
theorem :
\begin{theo}\label{main theo}
Let $\G$ be the Galois group of a Lie-irreducible generalized
$q$-hypergeometric equation $\hyp{\und a}{\und b}$ of order three
with $q$-real parameters. Then $\gzder=\sl{3}$. More precisely :
\begin{itemize}
  \item[$\bullet$] $\G=\gl{3}$ if $\frac{a_1a_2a_3}{b_2b_3} \not \in
  q^\rinteger$;
  \item[$\bullet$] $\G=\overline{\langle \sl{3}, e^{2\pi i (\beta_2+\beta_3)}, v_1v_2 \rangle}$ if $\frac{a_1a_2a_3}{b_2b_3} \in q^\rinteger$.
\end{itemize}
\end{theo}

As usually, $\gzder$ denotes the derived group of $\gz$, the neutral
component of $\G$.

\section{Galois theory for regular singular $q$-difference
equations}\label{rappels galois}

Using analytic tools together with Tannakian duality, Sauloy
developed in \cite{sauloyqgaloisfuchs} a Galois theory for regular
singular $q$-difference systems. In this section, we shall first
recall the principal notions used in \cite{sauloyqgaloisfuchs},
mainly the Birkhoff matrix and the twisted Birkhoff matrix. Then
we shall explain briefly that this lead to a Galois theory for
regular singular $q$-difference systems. Last, we shall state a
density theorem for these Galois groups, which will be of main
importance in the rest of the paper.

\subsection{Basic notions}\label{section the basic objects}

Let us consider $A \in Gl_n(\complex(\{z\}))$. Following Sauloy in
\cite{sauloyqgaloisfuchs} (see also
\cite{sauloysystqdiffsinguregu}), the $q$-difference system :
\begin{equation}\label{general q diff system}Y(qz)=A(z)Y(z)\end{equation}
is said to be \textit{Fuchsian} at $0$ if $A$ is holomorphic at
$0$ and if $A(0)\in \gl{n}$. Such a system is
\textit{non-resonant} at $0$ if, in addition, $Sp(A(0)) \cap
q^{\rinteger^*} Sp(A(0))=\emptyset$. Last we say that the above
$q$-difference system is \textit{regular singular} at $0$ if there
exists $R^{(0)}\in \gln (\complex(\{z\}))$ such that the
$q$-difference system defined by
$(R^{(0)}(qz))^{-1}A(z)R^{(0)}(z)$ is Fuchsian at $0$. We have
similar notions at $\infty$ using the variable change $z
\leftarrow 1/z$.

In case that the system is global, that is $A\in \gln
(\complex(z))$, we say that the system (\ref{general q diff
system}) is \textit{Fuchsian} (resp. \textit{Fuchsian and
non-resonnant}, \textit{regular singular}) if it is Fuchsian
(resp. Fuchsian and non-resonnant, regular singular) both at $0$
and at $\infty$.

For instance, the basic hypergeometric system (\ref{syst hypergeo}) is Fuchsian.\\

\textit{Local solution at $0$}. Suppose that (\ref{general q diff
system}) is Fuchsian and non-resonant at $0$ and consider
$J^{(0)}$ a Jordan normal form of $A(0)$. Thanks to
\cite{sauloyqgaloisfuchs} there exists $F^{(0)} \in \gln
(\complex\{z\})$ such that :
\begin{equation} \label{transfo jauge}
F^{(0)}(qz)J^{(0)}=A(z)F^{(0)}(z).
\end{equation}
Therefore, if $e^{(0)}_{J^{(0)}}$ denotes a fundamental system of
solutions of the $q$-difference system with constant coefficients
$X(qz)=J^{(0)}X(z)$, the matrix-valued function
$Y^{(0)}=F^{(0)}e^{(0)}_{J^{(0)}}$ is a fundamental system of
solutions of (\ref{general q diff system}). In
\cite{sauloyqgaloisfuchs}, $e^{(0)}_{J^{(0)}}$ is defined as
follows. We denote by $\theta_q$ the Jacobi theta function defined
by
$\theta_q(z)=\pochamer{q}{q}{\infty}\pochamer{z}{q}{\infty}\pochamer{q/z}{q}{\infty}$.
This is a meromorphic function over $\complex^*$ whose zeros are
simple and located on the discrete logarithmic spiral
$q^\rinteger$. Moreover, the functional equation
$\theta_q(qz)=-z^{-1}\theta_q(z)$ holds. Now we introduce, for all
$\lambda \in \complex^*$ such that $|q| \leq |\lambda| < 1$, the
$q$-character
$e^{(0)}_{\lambda}=\frac{\theta_q}{\theta_{q,\lambda}}$ with
$\theta_{q,\lambda}(z)=\theta_q(\lambda z)$ and we extend this
definition to an arbitrary nonzero complex number $\lambda \in
\complex^*$ requiring the equality
$e^{(0)}_{q\lambda}=ze^{(0)}_{\lambda}$. If
$D=P\text{diag}(\lambda_1,...,\lambda_n)P^{-1}$ is a semisimple
matrix then we set
$e^{(0)}_{D}:=P\text{diag}(e^{(0)}_{\lambda_1},...,e^{(0)}_{\lambda_n})P^{-1}$.
It is easily seen that it does not depend on the chosen
diagonalization. Furthermore, consider
$\ell_q(z)=-z\frac{\theta_q'(z)}{\theta_q(z)}$ and, if $U$ is a
unipotent matrix, $e^{(0)}_{U}=\Sigma_{k=0}^n \ell_q^{(k)}
(U-I_n)^k$ with $\ell_q^{(k)}=\binom{\ell_q}{k}$. Finally if
$J^{(0)}=D^{(0)}U^{(0)}$ is the multiplicative Dunford
decomposition of $J^{(0)}$, with $D^{(0)}$ semi-simple and
$U^{(0)}$ unipotent, we set
$e^{(0)}_{J^{(0)}}=e^{(0)}_{D^{(0)}}e^{(0)}_{U^{(0)}}$.\\

\textit{Local solution at $\infty$}. Using the variable change $z
\leftarrow 1/z$, we have a similar construction at $\infty$. The
corresponding fundamental system of solutions is denoted by
$Y^{(\infty)}=F^{(\infty)}e^{(\infty)}_{J^{(\infty)}}$.\\

Throughout this section we assume that (\ref{general q diff
system}) has coefficients in $\complex(z)$ and that it is Fuchsian
and non-resonant.\\

\textit{Birkhoff matrix}. The linear relations between the two
fundamental systems of solutions introduced above are given by the
Birkhoff matrix (also called connection matrix)
$P=(Y^{(\infty)})^{-1}Y^{(0)}$. Its entries are elliptic functions
\textit{i.e.} meromorphic functions over the elliptic curve
$\eq=\complex^* / q^\rinteger$. \\

\textit{Twisted Birkhoff matrix}. In order to describe a
Zariki-dense set of generators of the Galois group associated to
the system (\ref{general q diff system}), we introduce a
``twisted" connection matrix. According to
\cite{sauloyqgaloisfuchs}, we choose for all $z \in \complex^*$ a
group endomorphism $g_z$ of $\complex^*$ sending $q$ over $z$.
Before giving an explicit example, we have to introduce some
notations. Let us consider $\tau \in \complex$ such that
$q=e^{-2\pi i \tau}$ and set, for all $y \in \complex$,
$q^y=e^{-2\pi i \tau y}$. We also define the (non continuous)
function $\log_q$ on the whole punctured complex plane
$\complex^*$ by $\log_q(q^y)=y$ if $y\in \complex^* \setminus
q^\real$ and we require that its discontinuity is just before the
cut when turning counterclockwise around $0$. We can now give an
explicit example of endomorphism $g_z$ namely the function $g_z:
\complex^*=\mathbb{U} \times q^\real \rightarrow \complex^*$
sending $uq^\omega$ to $g_z(uq^\omega)=z^{\omega}=e^{-2\pi i \tau
\log_q(z) \omega}$ for $(u,\omega)\in \mathbb{U} \times \real$,
where $\mathbb{U}\subset \complex$ is the unit circle.

Then we set, for all $z$ in $\complex^*$,
$\psi_z^{(0)}(\lambda)=\frac{\qcar{\lambda}(z)}{g_z(\lambda)}$ and
we define  $\psi_z^{(0)}\left(D^{(0)}\right)$, the \textit{twisted
factor} at $0$, by
$\psi_z^{(0)}\left(D^{(0)}\right)=P\text{diag}(\psi_z^{(0)}(\lambda_1),...,\psi_z^{(0)}(\lambda_n))P^{-1}$
with $D^{(0)}=P\text{diag}(\lambda_1,...,\lambda_n)P^{-1}$. We
have a similar construction at $\infty$ by using the variable
change $z \leftarrow 1/z$. The corresponding twisting factor is
denoted by $\psi_z^{(\infty)}(J^{(\infty)})$.

Finally, the twisted connection matrix $\breve{P}(z)$ is :
\begin{eqnarray*}
\breve{P}(z)&=&\psi_z^{(\infty)}\left(D^{(\infty)}\right)P(z)\psi_z^{(0)}\left(D^{(0)}\right)^{-1}.
\end{eqnarray*}

\subsection{Definition of the difference Galois groups}\label{constr}

The definition of the Galois groups of regular singular
$q$-difference systems given by Sauloy in
\cite{sauloyqgaloisfuchs} is somewhat technical and long. Here we
do no more than describe the underlying idea.\\

\textit{(Global) Galois group.} Let us denote by $(\mathcal{E},\otimes)$ the rigid abelian tensor category
of regular singular $q$-difference systems with
coefficients in $\complex(z)$. This is actually a Tannakian category over $\complex$ but the existence of a fiber functor
is not obvious. We shall now explain how one can get a complex valued fiber functor \textit{via} a Riemann-Hilbert correspondance.

The category $\mathcal{E}$ is equivalent to
the category $\mathcal{C}$ of connection triples whose objects are
triples $(A^{(0)},P,A^{(\infty)}) \in \gln (\complex) \times \gln
(\mathcal{M}(\eq)) \times \gln (\complex)$ (we refer to
\cite{sauloyqgaloisfuchs} for the complete definition of
$\mathcal{C}$). Furthermore $\mathcal{C}$ can be endowed with a
tensor product $\und \otimes$ making the above equivalence of
categories compatible with the tensor products. Let us emphasize
that $\und \otimes$ is not the usual tensor product for matrices.
Indeed some twisting factors appear because of the bad
multiplicative properties of the $q$-characters $e_{q,c}$ : in
general $e_{q,c}e_{q,d}\neq e_{q,cd}$.

The category rigid abelian tensor category
$\mathcal{C}$ is a Tannakian category over $\complex$. The functor
$\omega_0$ from $\mathcal{C}$ to $Vect_\complex$ sending an object
$(A^{(0)},P,A^{(\infty)})$ to the underlying vector space
$\complex^n$ on which $A ^{(0)}$ acts is a fiber functor. Let us
remark that there is a similar fiber functor $\omega_\infty$ at
$\infty$. Following the general formalism of the theory of
Tannakian categories (see \cite{deligne}), the \textit{absolute
Galois group} of $\mathcal{C}$ (or, using the above equivalence of
categories, of $\mathcal E$) is defined as the pro-algebraic group
$\textsf{Aut}^{\und \otimes}(\omega_0)$ and the \textit{global
Galois group of an object $\chi$ }of $\mathcal{C}$ (or, using the
above equivalence of categories, of an object of $\mathcal{E}$) is
the complex linear algebraic group $\textsf{Aut}^{\und
\otimes}(\omega_{0|\langle \chi \rangle})$ where $\langle \chi
\rangle$ denotes the Tannakian subcategory of $\mathcal{C}$
generated by $\chi$. For the sake of simplicity, we will often
call $\textsf{Aut}^{\und \otimes}(\omega_{0|\langle \chi
\rangle})$ the \textit{Galois group} of $\chi$
(or, using the above equivalence of categories, of the corresponding object of $\mathcal{E}$). \\

\textit{Local Galois groups.} Notions of
local Galois groups at $0$ and at $\infty$ are also available. They
are subgroups of the (global) Galois group. Nevertheless, since
these groups are of second importance in what follows, we omit the
details and
we refer the interesting reader to \cite{sauloyqgaloisfuchs}.\\

\noindent \textbf{Remark.} In \cite{psgaloistheory}, Van der Put
and Singer showed that the Galois groups defined using a
Picard-Vessiot theory can be recover by means of Tannakian duality
: it is the group of tensor automorphisms of some complex
valued fiber functor over $\mathcal E$. Since two complex valued
fiber functors on a same Tannakian category are
isomorphic, we deduce that
Sauloy's and Van der Put and Singer's theories coincide.\\

In the rest of this section we exhibit some natural elements of
the Galois group of a given Fuchsian $q$-difference system and we
state a density theorem.

\subsection{The density theorem}

In what follows, we assume that the $q$-difference system
(\ref{general q diff system}) has coefficients in $\complex(z)$
and that it is Fuchsian and non-resonant.

Fix a ``base point" $y_0\in \bon=\complex^* \setminus
\{\text{zeros of } \text{det}(P(z)) \text{ or poles of } P(z)\}$ .
Sauloy exhibits in \cite{sauloyqgaloisfuchs} the following
elements of the (global) Galois group associated to the
$q$-difference system (\ref{general q diff system}) :
\begin{itemize}
\item[1.a.] $\gamma_1(D^{(0)})$ and $\gamma_2(D^{(0)})$ where :
$$\gamma_1:\complex^*=\mathbb{U} \times q^\real \rightarrow \mathbb{U}$$
is the projection over the first factor and : $$\gamma_2 :
\complex^*=\mathbb{U} \times q^\real \rightarrow \complex^*$$ is
defined by $\gamma_2(uq^\omega)=e^{2\pi i \omega}$.
\item[1.b.] $U^{(0)}$.
\item[2.a.] $\breve{P}(y_0)^{-1}\gamma_1(D^{(\infty)})\breve{P}(y_0)$ and $\breve{P}(y_0)^{-1}\gamma_2(D^{(\infty)})\breve{P}(y_0)$.
\item[2.b.] $\breve{P}(y_0)^{-1} U^{(\infty)} \breve{P}(y_0)$.
\item[3.] $\breve{P}(y_0)^{-1}\breve{P}(z)$, $z \in \bon$. \\
\end{itemize}

The following result is due to Sauloy \cite{sauloyqgaloisfuchs}.

\begin{theo}\label{dens theo}
The algebraic group generated by the matrices 1.a. to 3. is the
(global) Galois group $\G$ of the $q$-difference system
(\ref{general q diff system}). The algebraic group generated by
the matrices 1.a. and 1.b. is the local Galois group at $0$ of the
$q$-difference system (\ref{general q diff system}). The algebraic
group generated by the matrices 2.a. and 2.b. is the local Galois
group at $\infty$ of the $q$-difference system (\ref{general q
diff system}).
\end{theo}

The algebraic group generated by the matrices involved in 3. is
called the \textit{connection component} of the Galois group $\G$.


\section{Irreducibility and Lie-irreducibility}\label{lie irr}

Let us introduce some terminologies.

Let $\mathcal S$ be an object of $\mathcal E$. The fiber functor
$\omega_0$ induces a equivalence of tensor categories :
$$\langle \mathcal S \rangle \xrightarrow[]{\sim} (\text{finite dimensional }\complex\text{-representations of }\G=\textsf{Aut}^{\und \otimes}(\omega_{0|\langle \mathcal S \rangle})).$$

A regular singular $q$-difference system $\mathcal S$ of order $n$
is \textit{irreducible} if it corresponds to an irreducible
representation of its Galois group $\G$; it is
\textit{Lie-irreducible} if the restriction to $\gz$ of the
representation of $\G$ corresponding to $\mathcal S$ is irreducible.

We have the following obvious implication :
 $$\mathcal S\text{ Lie-irr. } \Longrightarrow \mathcal
S\text{ irr. }.$$

To any $q$-difference equation corresponds a $q$-difference system
by the usual trick (the converse is also true : this result is
known as the cyclic vector lemma), so that we can speak of the
Galois group of a $q$-difference equation. The irreducibility of a
$q$-difference operator $L$ in the above ``Galoisian" sense is
equivalent to the irreducibility of $L$
as an operator (immediate consequence of the Tannakian duality).

It is convenient and useful to introduce the notion
of \textit{$q$-difference module} (see
\cite{psgaloistheory,sauloynewton}). A $q$-difference module is a
module over the non-commutative algebra (quantum algebra)
$\mathcal D_q=\complex(z)\langle \sigq,\sigq^{-1} \rangle$ of
non-commutative polynomials satisfying the relation $\sigq z = qz
\sigq$. The $q$-difference module associated to the $q$-difference
operator $L\in\mathcal D_q$ is $\mathcal D_q / \mathcal D_q L$.
The irreducibility of $L$ as an operator is equivalent to the
simplicity of the corresponding $q$-difference module.

The counterpart of the following result for generalized
hypergeometric differential equations is due to Katz (see
\cite{katzexpsum}). The $q$-difference case is treated in
(\cite{}).

\begin{prop}\label{qhyp irr}
The generalized $q$-hypergeometric system $\hyp{\und a}{\und b}$
with parameters $(\und a,\und b) \in (\complex^*)^r \times
(\complex^*)^s$ is irreducible if and only if for all $i,j$ we have $a_i/b_j
\not \in q^\rinteger$. In this case, for all $(\und \alpha,\und
\beta) \in  \rinteger^r \times \rinteger^s$ the $q$-difference
modules associated to $\hyp{\und a}{\und b}$ and to
$\hyp{q^{\alpha_1}a_1, \cdots
,q^{\alpha_r}a_r}{q^{\beta_1}b_1,\cdots,q^{\beta_s}b_s}$ are
isomorphic (hence, they have isomorphic Galois groups).
\end{prop}

\section{Preliminaries on algebraic group theory}\label{pre alg gps}

Recall that $\sl{2}$ admits up to isomorphism exactly one
irreducible complex linear representation in each dimensions. The
unique (up to isomorphism) irreducible representation of degree
$3$ of $\sl{2}$ is $\textsf{Sym}^{\otimes 2}(\textsf{Std})$, the
symmetric square of the standard representation of $\sl{2}$ of
degree $2$, and is explicitly given by the following formula :
\begin{eqnarray}\label{rep sl}
\rho : \sl{2} & \rightarrow & \gl{3}\\
\nonumber \pmatrice{a}{b}{c}{d} & \mapsto &
\pmatricet{a^2}{2ab}{b^2}{ac}{ad+bc}{bd}{c^2}{2cd}{d^2}.
\end{eqnarray}
The image of this representation, denoted $\psl$, is an algebraic
subgroup of $\gl 3$ isomorphic to $\sl{2}/\{\pm I_2\}$. Of course,
this notation is not usual but will be very convenient.

For the following result, we refer to \cite{singulm} and to the references therein.

\begin{prop}\label{prop alter}
Let $\G$ be a connected semi-simple algebraic subgroup of $\sl{3}$
acting irreducibly on $\complex^3$. Then we have the following
alternative :
\begin{itemize}
\item[-] either $\G=\sl{3}$;
\item[-] or $\G$ is conjugate to $\psl$.
\end{itemize}
\end{prop}

The previous Proposition can be also seen as a special case of the
prime recognition lemma due to O. Gabber (see \cite{katzexpsum}).

Let us now list some properties of $\psl$
which will be useful in the rest of the paper. We set $\pgl=\complex^*
\cdot \psl$.

\begin{lem}\label{norm psl}
The normalizer of $\psl$ in $\gl{3}$ is $\pgl$.
\end{lem}
\begin{proof}
See \cite{singulm,singersolvhom}.
\end{proof}

\begin{lem}\label{lem diag psl} Let us consider $M \in \psl \setminus \{I_3\}$.

\begin{itemize}
  \item [(i)] The set of eigenvalues of $M$ is of the form
$\{1,\alpha,\alpha^{-1}\}$ with $\alpha \in \complex^*$;
  \item [(ii)] $M$ is
semi-simple iff $\alpha \neq 1$; in this case, $M$ is conjugated
to a diagonal matrix by an element of  $\psl$;
  \item [(iii)] if $M$ is not semi-simple (that is $\alpha=1$) then $M$ is
conjugated to \scriptsize $\pmatricet{1}{2}{1}{0}{1}{1}{0}{0}{1}$
\normalsize by an element of  $\psl$. 
\end{itemize}

\end{lem}
\begin{proof}
Let $N$ be an element of $\sl{2}$ such that $M=\rho(N)$. If $N$ is
semi-simple, it is conjugated to some
$\text{diag}(\lambda,\lambda^{-1})$. Hence $M=\rho(N)$ is
conjugated to $\text{diag}(\lambda^2,1,\lambda^{-2})$ by an
element of $\psl$. On the other hand, if $N$ is not semi-simple,
then $N$ is conjugated to \scriptsize $\pm \pmatrice{1}{1}{0}{1}$
\normalsize. Consequently, $M=\rho(N)$ is conjugated to
$\rho$\scriptsize $\left(\pm \pmatrice{1}{1}{0}{1}\right)$
$=\pmatricet{1}{2}{1}{0}{1}{1}{0}{0}{1}$\normalsize by an element
of $\psl$. The Lemma follows easily.
\end{proof}

We will also need the following elementary results, the proofs of
which are left to the reader. We denote by $\perm$ the finite
subgroup of $\gl{3}$ of permutation matrices; $\permc$ denotes the
subgroup of $\gl{3}$ generated by $\perm$ and by the invertible
diagonal matrices.

\begin{lem}\label{lem1}
Let us consider $M\in \gl{3}$ with three distinct eigenvalues.
Suppose that $P,P'\in \gl{3}$ are such that $P^{-1}MP$ and
$P'^{-1}MP'$ are diagonal, then there exists $\Sigma\in\permc$
such that $P'=P\Sigma$.
\end{lem}

\begin{lem}\label{lem2}
The normalizer of $\pmatricet{1}{1}{0}{0}{1}{1}{0}{0}{1}$ in
$\gl{3}$ is made of upper-triangular matrices.
\end{lem}

\section{Galois groups of Lie-irreducible equations}\label{calcul}

In order to prove our main Theorem \ref{main theo} we study the
non-resonant Lie-irreducible equations in several steps
corresponding to various logarithmic and non-logarithmic cases.

\subsection{Non-resonant and non-logarithmic case}\label{section
generique}

\textit{\textbf{Notations.}} It will be convenient to use the
following notations. Let $\und u$ be a $3$-uple of
$(\complex^*)^3$.
\begin{itemize}
\item[-] for all $\lambda \in \complex^*$, $\lambda \und u =(\lambda u_1,\lambda u_2,\lambda u_3)$ and $\lambda/\und
u=(\lambda/u_1,\lambda/u_2,\lambda/u_3)$;
\item[-]
$V(\und
u)=\pmatricet{1}{1}{1}{q/u_1}{q/u_2}{q/u_3}{(q/u_1)^2}{(q/u_2)^2}{(q/u_3)^2};$
  \item[-]
$\text{diag}(\und u)=\pmatricet{u_1}{0}{0}{0}{u_2}{0}{0}{0}{u_3};$
  \item[-] $\und u_{\check{1}}=(u_2,u_3)$, $\und
  u_{\check{2}}=(u_1,u_3)$ and $\und u_{\check{3}}=(u_1,u_2)$;
  \item[-] $\pochamer{\und u}{q}{\infty}=\lim_{n \rightarrow
  +\infty} \pochamer{\und u}{q}{n}$;
  \item[-] $a_i=u_iq^{\alpha_i}$ and $b_j=v_jq^{\beta_j}$ with
$u_i,v_j\in\mathbb{U}$ and $\alpha_i,\beta_j \in \real$ (we choose
a logarithm of $q$).
\end{itemize}

\noindent\textit{\textbf{Hypotheses.}} In this section we suppose
that : \begin{equation}\label{cond non res 1} \forall i\neq j,\ \
a_i/a_j \not \in q^\rinteger \text{ and } b_2/b_3,b_2,b_3 \not \in
q^\rinteger.\end{equation}

\noindent\textit{\textbf{Basic objects.}} Let us give the explicit
form of the matrices involved in the density Theorem \ref{dens theo}.\\

\textit{Local solution at $0$.} We have :
$$\qhypermatricetcour{\und{a}}{\und b}{0}=V(\und b) \jtz{\und b} V(\und b)^{-1}$$
with $\jtz{\und b}= \text{diag}(q/\und b)$.
Hence the system (\ref{syst hypergeo}) is non-resonant and
non-logarithmic at $0$ (that is $\qhypermatricetcour{\und a}{\und
b}{0}$ is semi-simple; hence we do not have to use the
$q$-logarithm $\ell_q$ for solving the $q$-difference system
defined by $\qhypermatricetcour{\und a}{\und b}{0}$). A
fundamental system of solutions at $0$ of (\ref{syst hypergeo}) as
described in section \ref{section the basic objects} is given by
$\ytz{\und a}{\und b}{z}=\ftz{\und a}{\und b}{z}e^{(0)}_{\jtz{\und
b}}(z)$ with :
$$\ftz{\und a}{\und b}{z}=\pmatricet{\qhypertcour{\und a}{\und b}{z}}
{\qhypertcour{\frac{q}{b_2} \und a}{\frac{q}{b_2} \und b _{}}{z}}
{\qhypertcour{\frac{q}{b_3} \und a}{\frac{q}{b_3} \und b _{}}{z}}
{\qhypertcour{\und a}{\und
b}{qz}}{\frac{q}{b_2}\qhypertcour{\frac{q}{b_2} \und
a}{\frac{q}{b_2} \und b _{}}{qz}}
{\frac{q}{b_3}\qhypertcour{\frac{q}{b_3} \und a}{\frac{q}{b_3}
\und b _{}}{qz}} {\qhypertcour{\und a}{\und
b}{q^2z}}{\left(\frac{q}{b_2}\right)^2\qhypertcour{\frac{q}{b_2}
\und a}{\frac{q}{b_2} \und b _{}}{q^2z}}
{\left(\frac{q}{b_3}\right)^2\qhypertcour{\frac{q}{b_3} \und a}{\frac{q}{b_3} \und b _{}}{q^2z}}.$$\\

\textit{Generators of the local Galois group at $0$.} We have two generators :
$$\text{diag}(e^{2 \pi i \und \beta}) \text{ and }
\text{diag}(\und v).$$


\textit{Local solution at $\infty$.} We have :
$$\qhypermatricetcour{\und a}{\und b}{\infty}=V(q\und a) \jtinf{\und a} V(q\und a)^{-1}$$
with $\jtinf{\und a}= \text{diag}(1/\und a)$.
Hence the system is non-resonant and non-logarithmic at $\infty$
and, a fundamental system of solutions at $\infty$ of (\ref{syst
hypergeo}) as described in section \ref{section the basic objects}
is given by $\ytinf{\und a}{\und b}{z}=\ftinf{\und a}{\und
b}{z}e^{(\infty)}_{\jtinf{\und a}}(z)$ with : \tiny $$\ftinf{\und
a}{\und
b}{z}=\pmatricet{\qhypertcour{a_1q/\und{b}}{a_1q/\und{a}_{}}{\frac{qb_2b_3}{a_1a_2a_3}z^{-1}}}
{\qhypertcour{a_2q/\und{b}}{a_2q/\und{a}_{}}{\frac{qb_2b_3}{a_1a_2a_3}z^{-1}}}
{\qhypertcour{a_3q/\und{b}}{a_3q/\und{a}_{}}{\frac{qb_2b_3}{a_1a_2a_3}z^{-1}}}
{\frac{1}{a_1}
\qhypertcour{a_1q/\und{b}}{a_1q/\und{a}_{}}{\frac{b_2b_3}{a_1a_2a_3}z^{-1}}}{\frac{1}{a_2}
\qhypertcour{a_2q/\und{b}}{a_2q/\und{a}_{}}{\frac{b_2b_3}{a_1a_2a_3}z^{-1}}}{\frac{1}{a_3}
\qhypertcour{a_3q/\und{b}}{a_3q/\und{a}_{}}{\frac{b_2b_3}{a_1a_2a_3}z^{-1}}}{\left(\frac{1}{a_1}\right)^2
\qhypertcour{a_1q/\und{b}}{a_1q/\und{a}_{}}{\frac{b_2b_3}{qa_1a_2a_3}z^{-1}}}{\left(\frac{1}{a_2}\right)^2
\qhypertcour{a_2q/\und{b}}{a_2q/\und{a}_{}}{\frac{b_2b_3}{qa_1a_2a_3}z^{-1}}}{\left(\frac{1}{a_3}\right)^2
\qhypertcour{a_3q/\und{b}}{a_3q/\und{a}_{}}{\frac{b_2b_3}{qa_1a_2a_3}z^{-1}}}
.$$\\
\normalsize

\textit{Generators of the local Galois group at $\infty$.} We have two generators :
$$\breve{P}(y_0)^{-1}\text{diag}(e^{2 \pi i \und \alpha})\breve{P}(y_0) \text{ and }
\breve{P}(y_0)^{-1}\text{diag}(\und u)\breve{P}(y_0).$$\\

\textit{Birkhoff matrix.} The Barnes-Mellin-Watson formula (see
\cite{gasperrahman} p. 110) entails that :

$$P(z)= (e^{(\infty)}_{\jtinf{\und a}}(z))^{-1}
\pmatricet {\frac{\pochamer{\und
a_{\check1},\und{b}_{\check1}/a_1}{q}{\infty}}{\pochamer{\und
b_{\check1},\und a_{\check1}/a_1}{q}{\infty}} \frac{\theta_q (a_1
z)}{\theta_q(z)}} {\frac{\pochamer{\frac{q}{b_2} \und
a_{\check1},\und
b_{\check2}/a_1}{q}{\infty}}{\pochamer{\frac{q}{b_2} \und
b_{\check2},\und a _{\check 1}/a_1}{q}{\infty}} \frac{\theta_q
(\frac{qa_1}{b_2} z)}{\theta_q(z)}} {\frac{\pochamer{\frac{q}{b_3}
\und a_{\check1},\und b _{\check
3}/a_1}{q}{\infty}}{\pochamer{\frac{q}{b_3}\und b_{\check3},\und a
_{\check 1}/a_1}{q}{\infty}} \frac{\theta_q (\frac{qa_1}{b_3}
z)}{\theta_q(z)}} {\frac{\pochamer{\und
a_{\check2},\und{b}_{\check1}/{a_2}}{q}{\infty}}{\pochamer{\und
b_{\check1},\und a_{\check2}/a_2}{q}{\infty}} \frac{\theta_q (a_2
z)}{\theta_q(z)}} {\frac{\pochamer{\frac{q}{b_2} \und
a_{\check2},\und b_{\check
2}/a_2}{q}{\infty}}{\pochamer{\frac{q}{b_2}\und b_{\check 2},\und
a _{\check 2}/a_2}{q}{\infty}} \frac{\theta_q (\frac{qa_2}{b_2}
z)}{\theta_q(z)}} {\frac{\pochamer{\frac{q}{b_3} \und
a_{\check2},\und b_{\check
3}/a_2}{q}{\infty}}{\pochamer{\frac{q}{b_3}\und b _{\check 3},\und
a _{\check 2}/a_2}{q}{\infty}} \frac{\theta_q (\frac{qa_2}{b_3}
z)}{\theta_q(z)}} {\frac{\pochamer{\und
a_{\check3},\und{b}_{\check1}/a_3}{q}{\infty}}{\pochamer{\und
b_{\check1},\und a_{\check3}/a_3}{q}{\infty}} \frac{\theta_q (a_3
z)}{\theta_q(z)}} {\frac{\pochamer{\frac{q}{b_2} \und
a_{\check3},\und
b_{\check2}/a_3}{q}{\infty}}{\pochamer{\frac{q}{b_2} \und
b_{\check2},\und a _{\check 3}/a_3}{q}{\infty}} \frac{\theta_q
(\frac{qa_3}{b_2} z)}{\theta_q(z)}} {\frac{\pochamer{\frac{q}{b_3}
\und a_{\check3},\und b_{\check3}/a_3}{q}{\infty}}
{\pochamer{\frac{q}{b_3}\und b_{\check3},\und a _{\check
3}/a_3}{q}{\infty}} \frac{\theta_q (\frac{qa_3}{b_3}
z)}{\theta_q(z)}}e^{(0)}_{\jtz{\und b}}(z).$$\\

\textit{Twisted Birkhoff matrix.} We deduce that :
$$\breve{P}(z)= \text{diag}((1/z)^{-\und \alpha})
\pmatricet {\frac{\pochamer{\und
a_{\check1},\und{b}_{\check1}/a_1}{q}{\infty}}{\pochamer{\und
b_{\check1},\und a_{\check1}/a_1}{q}{\infty}} \frac{\theta_q (a_1
z)}{\theta_q(z)}} {\frac{\pochamer{\frac{q}{b_2} \und
a_{\check1},\und
b_{\check2}/a_1}{q}{\infty}}{\pochamer{\frac{q}{b_2} \und
b_{\check2},\und a _{\check 1}/a_1}{q}{\infty}} \frac{\theta_q
(\frac{qa_1}{b_2} z)}{\theta_q(z)}} {\frac{\pochamer{\frac{q}{b_3}
\und a_{\check1},\und b _{\check
3}/a_1}{q}{\infty}}{\pochamer{\frac{q}{b_3}\und b_{\check3},\und a
_{\check 1}/a_1}{q}{\infty}} \frac{\theta_q (\frac{qa_1}{b_3}
z)}{\theta_q(z)}} {\frac{\pochamer{\und
a_{\check2},\und{b}_{\check1}/{a_2}}{q}{\infty}}{\pochamer{\und
b_{\check1},\und a_{\check2}/a_2}{q}{\infty}} \frac{\theta_q (a_2
z)}{\theta_q(z)}} {\frac{\pochamer{\frac{q}{b_2} \und
a_{\check2},\und b_{\check
2}/a_2}{q}{\infty}}{\pochamer{\frac{q}{b_2}\und b_{\check 2},\und
a _{\check 2}/a_2}{q}{\infty}} \frac{\theta_q (\frac{qa_2}{b_2}
z)}{\theta_q(z)}} {\frac{\pochamer{\frac{q}{b_3} \und
a_{\check2},\und b_{\check
3}/a_2}{q}{\infty}}{\pochamer{\frac{q}{b_3}\und b _{\check 3},\und
a _{\check 2}/a_2}{q}{\infty}} \frac{\theta_q (\frac{qa_2}{b_3}
z)}{\theta_q(z)}} {\frac{\pochamer{\und
a_{\check3},\und{b}_{\check1}/a_3}{q}{\infty}}{\pochamer{\und
b_{\check1},\und a_{\check3}/a_3}{q}{\infty}} \frac{\theta_q (a_3
z)}{\theta_q(z)}} {\frac{\pochamer{\frac{q}{b_2} \und
a_{\check3},\und
b_{\check2}/a_3}{q}{\infty}}{\pochamer{\frac{q}{b_2} \und
b_{\check2},\und a _{\check 3}/a_3}{q}{\infty}} \frac{\theta_q
(\frac{qa_3}{b_2} z)}{\theta_q(z)}} {\frac{\pochamer{\frac{q}{b_3}
\und a_{\check3},\und b_{\check3}/a_3}{q}{\infty}}
{\pochamer{\frac{q}{b_3}\und b_{\check3},\und a _{\check
3}/a_3}{q}{\infty}} \frac{\theta_q (\frac{qa_3}{b_3}
z)}{\theta_q(z)}} \text{diag}(z^{-\und \beta}).$$

$_{}$\vskip 10 pt

It will be convenient to write :
\begin{eqnarray*} p_{i,j}=p_{i,j}(\und a,\und b)= \frac{\pochamer{\frac{q}{b_j} \und
a_{\check i},\und{b}_{\check
j}/a_i}{q}{\infty}}{\pochamer{\frac{q}{b_j}\und b _{\check j},\und
a_{\check i}/a_i}{q}{\infty}}.
\end{eqnarray*}

For later use, we first compute the determinant and the minors of
$\breve{P}(z)$.

\begin{prop}\label{calcul det et mineurs}
We have :
\begin{itemize}
  \item[(i)] the minor $(i_1,i_2) \times (j_1,j_2)$ of $\breve{P}(z)$ is equal
  to :\small \begin{eqnarray}\nonumber \kappa_{(i_1,i_2) \times (j_1,j_2)} &=&
\frac{-q}{\pochamer{q}{q}{\infty}^2}\frac{a_{i_2}}{b_{j_1}}
\frac{\pochamer{\frac{q}{b_{j_1}}  a_{ i_3}, b_{
j_3}/a_{i_1},\frac{q}{b_{j_2}}  a_{ i_3}, b_{
j_3}/a_{i_2}}{q}{\infty}\theta_q (\frac{a_{i_1}}{a_{i_2}})\theta_q
(\frac{b_{j_1}}{b_{j_2}})}{\pochamer{\frac{q}{b_{j_1}}\und
b_{\check j_1},\und a _{\check i_1}/a_{i_1},\frac{q}{b_{j_2}}\und
b_{\check j_2},\und a _{\check
i_2}/a_{i_2}}{q}{\infty}} \cdot \\
\nonumber &&\ \ \ \ \ \ \ \ \ \ \ \ \ \ \ \ \ \ \ \ \ \ \ \ \ \ \ \ \ \ \ \ \ \ \ \cdot (1/z)^{-(\alpha_{i_1}+\alpha_{i_2})}z^{-(\beta_{j_1}+\beta_{j_2})}\frac{\theta_q
(\frac{q^2a_{i_1}a_{i_2}}{b_{j_1}b_{j_2}} z)}{\theta_q(z)}.\\
\end{eqnarray} \normalsize
  \item[(ii)] the determinant of $\breve{P}(z)$ is equal to :
\small  \begin{eqnarray}\label{det connec} \det(\breve{P}(z))&=&
-q\frac{(1-q/b_2)(1-q/b_3)(1/b_2-1/b_3)}{(1/a_2-1/a_1)(1/a_3-1/a_1)(1/a_2-1/a_3)}\cdot \nonumber \\
&&\ \ \ \ \ \ \ \ \ \ \ \ \ \ \ \ \ \ \ \ \ \ \ \ \cdot
(1/z)^{-(\alpha_1+\alpha_2+\alpha_3)}z^{-(\beta_1+\beta_2+\beta_3)}\frac{\theta_q(\frac{q^2a_1a_2a_3}{b_2b_3}
z)}{\theta_q(z)}.
\end{eqnarray} \normalsize
\end{itemize}
\end{prop}

\begin{proof}
\begin{itemize}
\item[(i)]The minor $(i_1,i_2) \times (j_1,j_2)$ of $\breve{P}(z)$ is equal
  to :
\small
\begin{eqnarray*}
&&\frac{(1/z)^{-(\alpha_{i_1}+\alpha_{i_2})}z^{-(\beta_{j_1}+\beta_{j_2})}}{\pochamer{\frac{q}{b_{j_1}}\und
b_{\check j_1},\und a _{\check i_1}/a_{i_1},\frac{q}{b_{j_2}}\und
b_{\check j_2},\und a _{\check i_2}/a_{i_2}}{q}{\infty}} \cdot \\
&& \ \ \ \ \ \ \ \ \ \ \ \cdot \left(\pochamer{\frac{q}{b_{j_1}}
\und a_{\check i_1},\und b_{\check j_1}/a_{i_1},\frac{q}{b_{j_2}}
\und a_{\check i_2},\und b_{\check
j_2}/a_{i_2}}{q}{\infty}\frac{\theta_q (\frac{qa_{i_1}}{b_{j_1}}
z)}{\theta_q(z)}\frac{\theta_q (\frac{qa_{i_2}}{b_{j_2}}
z)}{\theta_q(z)}\right. \\ &&\ \ \ \ \ \ \ \ \ \ \ \ \ \ \ \ \ \ \
\ \ \ \ \ \ \ \ \ \ \ \ \ \ \ \ \ \ \ \ \    - \left.
\pochamer{\frac{q}{b_{j_1}} \und a_{\check i_2},\und b_{\check
j_1}/a_{i_2},\frac{q}{b_{j_2}} \und a_{\check i_1},\und b_{\check
j_2}/a_{i_1}}{q}{\infty}\frac{\theta_q (\frac{qa_{i_2}}{b_{j_1}}
z)}{\theta_q(z)} \frac{\theta_q (\frac{qa_{i_1}}{b_{j_2}}
z)}{\theta_q(z)} \right).
\end{eqnarray*} \normalsize
Remark that the function : \small
\begin{eqnarray*} f(z)&=&\pochamer{\frac{q}{b_{j_1}}
\und a_{\check i_1},\und b_{\check j_1}/a_{i_1},\frac{q}{b_{j_2}}
\und a_{\check i_2},\und b_{\check
j_2}/a_{i_2}}{q}{\infty}\theta_q (\frac{qa_{i_1}}{b_{j_1}}
z)\theta_q (\frac{qa_{i_2}}{b_{j_2}} z)
\\ &&\ \ \ \ \ \ \ \ \ \ \ \ \ \ \ \ \ \ \ \ \ \ \ \ \ \ \ \ \ \ \
\ \ \ \ \     - \pochamer{\frac{q}{b_{j_1}} \und a_{\check
i_2},\und b_{\check j_1}/a_{i_2},\frac{q}{b_{j_2}} \und a_{\check
i_1},\und b_{\check j_2}/a_{i_1}}{q}{\infty}\theta_q
(\frac{qa_{i_2}}{b_{j_1}} z)\theta_q (\frac{qa_{i_1}}{b_{j_2}}
z)\end{eqnarray*} \normalsize vanishes for $z\in q^\rinteger$ and
$z\in \frac{b_{j_1}b_{j_2}}{q^2a_{i_1}a_{i_2}}q^\rinteger$.
Moreover, both functions $f/\theta_q^2$ and $\frac{\theta_q
(\frac{q^2a_{i_1}a_{i_2}}{b_{j_1}b_{j_2}} z)}{\theta_q(z)}$
satisfy the same linear homogeneous $q$-difference equation of
order one. We deduce that there exists a constant $\kappa \in
\complex^*$ such that, for all $z\in \complex^*$,
$(f/\theta_q^2)(z)=\kappa \frac{\theta_q
(\frac{q^2a_{i_1}a_{i_2}}{b_{j_1}b_{j_2}} z)}{\theta_q(z)}$. We
get the value of $\kappa$ by evaluating the above equality at
$z=\frac{b_{j_1}}{qa_{i_2}}$; we obtain :

\begin{eqnarray*}
\kappa&=&\frac{\pochamer{\frac{q}{b_{j_1}} \und a_{\check
i_1},\und b_{\check j_1}/a_{i_1},\frac{q}{b_{j_2}} \und a_{\check
i_2},\und b_{\check j_2}/a_{i_2}}{q}{\infty}\theta_q
(\frac{a_{i_1}}{a_{i_2}})\theta_q
(\frac{b_{j_1}}{b_{j_2}})}{\theta_q
(\frac{b_{j_1}}{qa_{i_2}})\theta_q
(q\frac{a_{i_1}}{b_{j_2}})}\\
&=&\frac{-q}{\pochamer{q}{q}{\infty}^2}\frac{a_{i_2}}{b_{j_1}}
\pochamer{\frac{q}{b_{j_1}}  a_{ i_3}, b_{
j_3}/a_{i_1},\frac{q}{b_{j_2}}  a_{ i_3}, b_{
j_3}/a_{i_2}}{q}{\infty}\theta_q (\frac{a_{i_1}}{a_{i_2}})\theta_q
(\frac{b_{j_1}}{b_{j_2}})\\
\end{eqnarray*}

Therefore the minor $(i_1,i_2) \times (j_1,j_2)$ of $\breve{P}(z)$
is equal to : \begin{eqnarray*} &&
\frac{-q}{\pochamer{q}{q}{\infty}^2}\frac{a_{i_2}}{b_{j_1}}
\frac{\pochamer{\frac{q}{b_{j_1}}  a_{ i_3}, b_{
j_3}/a_{i_1},\frac{q}{b_{j_2}}  a_{ i_3}, b_{
j_3}/a_{i_2}}{q}{\infty}\theta_q (\frac{a_{i_1}}{a_{i_2}})\theta_q
(\frac{b_{j_1}}{b_{j_2}})}{\pochamer{\frac{q}{b_{j_1}}\und
b_{\check j_1},\und a _{\check i_1}/a_{i_1},\frac{q}{b_{j_2}}\und
b_{\check j_2},\und a _{\check i_2}/a_{i_2}}{q}{\infty}} \cdot \\
&& \ \ \ \ \ \ \ \ \ \ \ \ \ \ \ \ \ \ \ \ \ \ \ \ \ \ \ \ \ \ \ \
\ \ \ \ \ \ \ \ \ \ \ \ \ \ \ \ \ \ \ \ \cdot
(1/z)^{-(\alpha_{i_1}+\alpha_{i_2})}z^{-(\beta_{j_1}+\beta_{j_2})}\frac{\theta_q
(\frac{q^2a_{i_1}a_{i_2}}{b_{j_1}b_{j_2}} z)}{\theta_q(z)}\\
\end{eqnarray*}

\item[(ii)] The following equality holds :
$\text{det}(\breve{P}(z))=(1/z)^{-(\alpha_1+\alpha_2+\alpha_3)}z^{-(\beta_1+\beta_2+\beta_3)}
\frac{f^{(0)}(\und a;\und b;z)}{f^{(\infty)}(\und a;\und b;z)}$
with $f^{(0)}(\und a;\und b;z)=\text{det}(\ftz{\und a}{\und
b}{z})$ and $f^{(\infty)}(\und a;\und b;z)=\text{det}(\ftinf{\und
a}{\und b}{z})$. From the functional equations verified by
$f^{(0)}(\und a;\und b;z)$ and $f^{(\infty)}(\und a;\und b;z)$,
and from the fact that these functions are respectively germs of
holomorphic functions at $0$ and at $\infty$, we deduce that
$\frac{f^{(0)}(\und a;\und b;z)}{f^{(\infty)}(\und a;\und b;z)}$
is a meromorphic function over $\complex^*$ whose poles are
simples and located on $q^\rinteger$ and whose zeros are also
simples and located on $\frac{b_2b_3}{q^2a_1a_3a_3}q^\rinteger$.
Furthermore we have : $\frac{f^{(0)}(\und a;\und
b;qz)}{f^{(\infty)}(\und a;\und
b;qz)}=\frac{b_2b_3}{q^2a_1a_3a_3}\frac{f^{(0)}(\und a;\und
b;z)}{f^{(\infty)}(\und a;\und b;z)}$. This ensures that there
exists $\eta \in\complex^*$ such that, for all $z\in\complex^*$,
$\frac{f^{(0)}(\und a;\und b;z)}{f^{(\infty)}(\und a;\und
b;z)}=\eta \frac{\theta_q(\frac{q^2a_1a_3a_3}{b_2b_3}
z)}{\theta_q(z)}$. We are going to determine the explicit value of
$\eta$.

Expanding $\text{det}(\breve{P}(z))$ with respect to the first
column, we get :
\small
\begin{eqnarray*}
\text{det}(\breve{P}(z))&=&
(1/z)^{-(\alpha_1+\alpha_2+\alpha_3)}z^{-(\beta_1+\beta_2+\beta_3)}\left(
\kappa_{(2,3)\times(2,3)} \frac{\theta_q
(\frac{q^2a_{2}a_{3}}{b_{2}b_{3}} z)}{\theta_q(z)}\frac{\theta_q
(a_1 z)}{\theta_q(z)}\right. \\
&& \ \ \ \ \ \ \ \ \ \ \ \ \ \ \ \ \ \ \ \ \ \ \ \ \ \ \ \ \ \ \ \ \ \ \ \ \ \ \ \ - \kappa_{(1,3)\times(2,3)}\frac{\theta_q
(\frac{q^2a_{1}a_{3}}{b_{2}b_{3}})}{\theta_q(z)} \frac{\theta_q
(a_2 z)}{\theta_q(z)}
\\ && \ \ \ \ \ \ \ \ \ \ \ \ \ \ \ \ \ \ \ \ \ \ \ \ \ \ \ \ \ \ \ \ \ \ \ \ \ \ \ \ \left. + \kappa_{(1,2)\times(2,3)} \frac{\theta_q
(\frac{q^2a_{1}a_{2}}{b_{2}b_{3}})}{\theta_q(z)} \frac{\theta_q
(a_3 z)}{\theta_q(z)}\right).
\end{eqnarray*}
\normalsize
 We
obtain the value of $\eta$ by specializing the above equality at
$z=1/a_3$. Indeed a straightforward calculation shows that :
\small
\begin{eqnarray*} \eta
 \frac{\theta_q (\frac{q^2a_1a_2}{b_2b_3})}{\theta_q(1/a_3)} &=&
 \frac{q}{\pochamer{q}{q}{\infty}^6} \frac{1}{\pochamer{\und b_{\check1},\und
a_{\check1}/a_1,\frac{q}{b_2}\und b_{\check 2},\und a _{\check
2}/a_2,\frac{q}{b_3}\und b_{\check3},\und a _{\check
3}/a_3}{q}{\infty}}\frac{a_2^2a_3}{b_2}
\frac{\theta_q(\frac{b_2}{b_3})\theta_q(\frac{a_3}{a_2})\theta_q(\frac{a_1}{a_3})}{\theta_q(a_3)}
 \\
&& \ \ \ \ \ \ \ \ \ \ \ \ \ \ \ \ \ \ \ \ \ \ \ \ \cdot
\left(\theta_q(a_2)\theta_q(\frac{b_2}{a_1})\theta_q(\frac{b_3}{a_1})\theta_q(\frac{q^2a_2}{b_2b_3})\right.
\left.
-\theta_q(a_1)\theta_q(\frac{b_2}{a_2})\theta_q(\frac{b_3}{a_2})\theta_q(\frac{q^2a_1}{b_2b_3})\right).
\end{eqnarray*} \normalsize
Moreover : \small
$$\theta_q(a_2)\theta_q(\frac{b_2}{a_1})\theta_q(\frac{b_3}{a_1})\theta_q(\frac{q^2a_2}{b_2b_3})-
\theta_q(a_1)\theta_q(\frac{b_2}{a_2})\theta_q(\frac{b_3}{a_2})\theta_q(\frac{q^2a_1}{b_2b_3})=-a_2
\theta_q(b_3)\theta_q(\frac{q^2a_1a_2}{b_2b_3})\theta_q(\frac{a_1}{a_2})\theta_q(b_2)
$$ \normalsize (use the above formula for $f$ with a suitable choice of parameters).
Hence, we have the formula : \small \begin{eqnarray*}
\text{det}(\breve{P}(z))&=&
-q\frac{(1-q/b_2)(1-q/b_3)(1/b_2-1/b_3)}{(1/a_2-1/a_1)(1/a_3-1/a_1)(1/a_2-1/a_3)} \nonumber \\
&&\ \ \ \ \ \ \ \ \ \ \ \ \ \ \ \ \ \ \ \ \ \ \ \ \cdot
(1/z)^{-(\alpha_1+\alpha_2+\alpha_3)}z^{-(\beta_1+\beta_2+\beta_3)}\frac{\theta_q(\frac{q^2a_1a_2a_3}{b_2b_3}
z)}{\theta_q(z)}.
\end{eqnarray*} \normalsize

\end{itemize}

\end{proof}

\begin{theo}\label{lie irred prem}
Suppose that the $q$-hypergeometric equation $\hyp{\und a}{\und
b}$ is Lie-irreducible and has $q$-real
parameters. Then $\gzder=\sl{3}$. Moreover, we have :
\begin{itemize}
  \item[$\bullet$] $\G=\gl{3}$ if $\frac{a_1a_2a_3}{b_2b_3} \not \in
  q^\rinteger$;
  \item[$\bullet$] $\G=\overline{\langle \sl{3}, e^{2\pi i (\beta_2+\beta_3)}, v_1v_2 \rangle}$ if $\frac{a_1a_2a_3}{b_2b_3} \in q^\rinteger$.
\end{itemize}
\end{theo}

\begin{proof}
Since $\gz$ acts irreducibly and faithfully on $\complex^3$, $\gz$
is reductive. The general theory of algebraic groups entails that
$\gz$ is generated by its center $Z(\gz)$ together with its
derived subgroup $\gzder$ which is semi-simple and that $Z(\gz)$
acts as scalars. Hence, the connected semi-simple algebraic group
$\gzder \subset \sl{3}$ also acts irreducibly 
on
$\complex^3$. Therefore, Proposition \ref{prop alter} ensures that
$\gzder$ is either conjugated to $\psl$ or equal to $\sl{3}$.

Suppose that $\gzder$ is conjugated to $\psl$. Then, the Galois
group $\G$ being a subgroup of the normalizer of $\gzder$, we
deduce from Lemma \ref{norm psl} that $\G$ is a subgroup of some
conjugate of $\pgl$. Let $R \in \gl{3}$ be such that $\G \subset
R^{-1} \pgl R$.

The non-resonance hypothesis implies that $\sharp \{ e^{2\pi i
\beta_1}, e^{2\pi i \beta_2}, e^{2\pi i \beta_3}\}=3$ (resp.
$\sharp \{ e^{2\pi i \alpha_1}, e^{2\pi i \alpha_2}, e^{2\pi i
\alpha_3}\}=3$), that is, that the matrix $R\text{diag}(e^{2 \pi i
\und\beta})R^{-1}$ (resp. $R\breve{P}(z)^{-1}\text{diag}(e^{2 \pi
i \und\alpha})\breve{P}(z)R^{-1}$ for all $z\in\bon$) has three
distinct eigenvalues and hence, thanks to Lemma \ref{lem diag
psl}, is conjugated to some diagonal matrix (with three distinct
diagonal entries) by some $\widetilde{R}\in \psl$ (resp.
$\widetilde{R}_z\in \psl$). It follows from Lemma \ref{lem1} that
there exists $\Sigma \in \permc$ (resp. $\Sigma_z \in \permc$)
such that $R=\widetilde{R}\Sigma$ (resp.
$R\breve{P}(z)^{-1}=\widetilde{R}_z \Sigma_z$).


Consequently, for all $z\in\bon$, $\Sigma_z\breve{P}(z)\Sigma^{-1}
\in \psl$. Let us consider a set $\bon'\subset \bon$ with at least
one accumulation point in $\complex^*\setminus q^\rinteger$ such
that there exists $\Sigma'\in\perm$ such that, for all
$z\in\bon'$, there exists
$(\lambda_{1,z},\lambda_{2,z},\lambda_{3,z})\in (\complex^*)^3$
with
$\Sigma_z=\text{diag}(\lambda_{1,z},\lambda_{2,z},\lambda_{3,z})\Sigma'$.
In order to simplify the proof, we assume that
$\Sigma=\Sigma'=I_3$ : the reader will easily adapt the rest of
the proof to the general case (permutation of the indices). In
view of (\ref{rep sl}), we get that the entries of a matrix
$M=(m_{i,j})\in \psl$ satisfy $m_{1,2}^2=4m_{1,1}m_{1,3}$.
Applying this formula to $M=\Sigma_z\breve{P}(z)\Sigma^{-1}$, we
get a relation of the form :\scriptsize
\begin{eqnarray}\label{eq fonc}\left(\frac{\pochamer{\frac{q}{b_2} \und a_{\check1},\und
b_{\check2}/a_1}{q}{\infty}}{\pochamer{\frac{q}{b_2} \und
b_{\check2},\und a _{\check 1}/a_1}{q}{\infty}} \frac{\theta_q
(\frac{qa_1}{b_2}
z)}{\theta_q(z)}(1/z)^{-\alpha_1}z^{-\beta_2}\right)^2=cst \cdot
\frac{\pochamer{\und
a_{\check1},\und{b}_{\check1}/a_1,\frac{q}{b_3} \und
a_{\check1},\und b _{\check 3}/a_1}{q}{\infty}}{\pochamer{\und
b_{\check1},\und a_{\check1}/a_1,\frac{q}{b_3}\und
b_{\check3},\und a _{\check 1}/a_1}{q}{\infty}} \frac{\theta_q
(a_1 z)}{\theta_q(z)} \frac{\theta_q (\frac{qa_1}{b_3}
z)}{\theta_q(z)}(1/z)^{-2\alpha_1}z^{-(\beta_1+\beta_3)},\end{eqnarray}\normalsize
for all $z\in\bon'$, where $cst\in\complex^*$ does not depend on
$z$. Since $\bon'$ has an accumulation point in $\complex^*
\setminus q^\rinteger$, the principle of analytic continuation
entails that (\ref{eq fonc}) holds for all
$z\in\widetilde{\complex^*}$. Note that the Pochhammer symbols
involved in the above functional equation are nonzero (if not,
there would exist $i,j$ such that $a_i/b_j\in q^\rinteger$ : this
is excluded since in this case the system is reducible). Now, the
localization of the zeros of both sides of (\ref{eq fonc}) leads
us to a contradiction. Indeed, the left hand side vanishes on
$\frac{qa_1}{b_2}q^\rinteger$ whereas the right hand side vanishes
exactly on $a_1q^\rinteger$ and on $\frac{qa_1}{b_3}q^\rinteger$ :
this is a contradiction since the non-resonance hypothesis
(\ref{cond non res 1}) implies that these three discrete
$q$-logarithmic spirals are distinct.

Therefore $\gzder = \sl{3}$. The theorem follows from this and
from formula (\ref{det connec}).
\end{proof}

\subsection{$b_2=b_3 \not \in q^\rinteger$ and the system is non-resonant and non-logarithmic at $\infty$}

\noindent\textit{\textbf{Hypotheses.}} In this section we suppose
that : \begin{equation}\label{cond non res 2}\forall i\neq j,\ \
a_i/a_j \not \in q^\rinteger \text{ and } b_2=b_3 \not \in
q^\rinteger.\end{equation}

\noindent\textit{\textbf{Basic objects.}}

It is easily seen that $\qhypermatricetcour{\und{a}}{\und b}{0}$
is conjugated to $\jtz{\und b}=
\pmatricet{1}{0}{0}{0}{q/b_2}{1}{0}{0}{q/b_2}$.\\

\textit{Generators of the local Galois group at $0$.} We have
three generators :
$$\text{diag}(e^{2 \pi i \und \beta}),~
\text{diag}(\und v) \text{ and }
\pmatricet{1}{0}{0}{0}{1}{b_2/q}{0}{0}{1}.$$

In this section we do not need neither the explicit form of the
local component at $\infty$ nor that of the twisted Birkhoff
matrix in order to determine the derived group of the neutral
component of the Galois group of a Lie-irreducible equation.

\begin{theo}
Let us suppose that the $q$-hypergeometric equation $\hyp{\und
a}{\und b}$ is Lie-irreducible. Then $\gzder=\sl{3}$.
\end{theo}

\begin{proof}
Arguing as for the proof of Theorem \ref{lie irred prem}, we get
that $\gzder$ is either conjugate to $\psl$, or equal to $\sl{3}$
and that, in case that $\gzder$ is conjugated to $\psl$, the
Galois group $\G$ is a subgroup of some conjugate of $\pgl$. Such
an inclusion is impossible since $\pgl$ contains exclusively
semi-simple matrices or matrices with exactly one Jordan bloc, so
that $\pgl$ cannot contain a conjugate of \scriptsize
$\pmatricet{1}{0}{0}{0}{1}{b_2/q}{0}{0}{1}$\normalsize.
\end{proof}

\subsection{$b_2=b_3=q$ and the system is non-resonant and non-logarithmic at $\infty$}

\noindent\textit{\textbf{Notations.}} We set $\und
q=(q,q,q)\in(\complex^*)^3$.\\

\noindent\textit{\textbf{Hypotheses.}} In this section we suppose
that : \begin{equation}\label{cond non res 3}\forall i\neq j,\ \
a_i/a_j \not \in q^\rinteger \text{ and } b_2=b_3=q.\end{equation}

\noindent\textit{\textbf{Basic objects.}}

\textit{Local solution at $0$.} We have :
$$\qhypermatricetcour{\und a}{\und q}{0}= \pmatricet{1}{-1}{1}{1}{0}{0}{1}{1}{0} \pmatricet{1}{1}{0}{0}{1}{1}{0}{0}{1} \pmatricet{1}{-1}{1}{1}{0}{0}{1}{1}{0}^{-1}.$$
Consequently, we are in the non-resonant logarithmic case at $0$.
In order to get a fundamental system of solutions at $0$ of
(\ref{syst hypergeo}) when $\und b=\und q$, we use a degeneration
procedure.

We consider the limit as $\und b$ tends to $\und q$ of the
following matrix valued function :

\begin{eqnarray*}
&&\ftz{\und a}{\und b}{z}V(\und b)^{-1}\pmatricet{1}{-1}{1}{1}{0}{0}{1}{1}{0}\\
&=& \ftz{\und a}{\und b}{z}
\pmatricet{1}{\frac{1-q^2/(b_2b_3)}{(\frac{q}{b_2}-1)(\frac{q}{b_3}-1)}}{\frac{q^2/(b_2b_3)}{(\frac{q}{b_2}-1)(\frac{q}{b_3}-1)}}
{0}{\frac{\frac{q}{b_3}-1}{(\frac{q}{b_3}-\frac{q}{b_2})(\frac{q}{b_2}-1)}}{\frac{-\frac{q}{b_3}}{(\frac{q}{b_3}-\frac{q}{b_2})(\frac{q}{b_2}-1)}}
{0}{\frac{\frac{q}{b_2}-1}{(\frac{q}{b_2}-\frac{q}{b_3})(\frac{q}{b_3}-1)}}{\frac{-\frac{q}{b_2}}{(\frac{q}{b_2}-\frac{q}{b_3})(\frac{q}{b_3}-1)}}.\\
\end{eqnarray*}
An easy computation shows that this limit does exist and we denote
it by $\ftz{\und a}{\und q}{z}$. From (\ref{transfo jauge}) we
deduce that $\ftz{\und a}{\und q}{z}$ satisfies $\ftz{\und a}{\und
q}{qz}\jtz{\und q}=\qhypermatricetcour{\und a}{\und q}{z}\ftz{\und
a}{\und q}{z}$ with $\jtz{\und
q}=\pmatricet{1}{1}{0}{0}{1}{1}{0}{0}{1}$. Hence, the matrix
$\ftz{\und a}{\und q}{z}$ being invertible as a matrix in the
field of meromorphic functions, the matrix valued function
$\ytz{\und a}{\und q}{z}=\ftz{\und a}{\und q}{z}e^{(0)}_{\jtz{\und
q}}(z)$ is a fundamental system of solutions of the
$q$-hypergeometric equation with $\und b=\und q$. Recall that $e^{(0)}_{\jtz{\und q}}(z)= \pmatricet{1}{\ell_q(z)}{\frac{\ell_q(z)(\ell_q(z)-1)}{2}}{0}{1}{\ell_q(z)}{0}{0}{1}$. \\
\\

\textit{Generators of the local Galois group at $0$.} We have the following generator :
$$\pmatricet{1}{1}{0}{0}{1}{1}{0}{0}{1}.$$\\

\textit{Local solution at $\infty$.} The situation is the same
than in section \ref{section generique}. We are in the
non-resonant and non-logarithmic case at $\infty$ and a
fundamental system of solutions at $\infty$ of (\ref{syst
hypergeo}) as described in section \ref{section the basic objects}
is given by $\ytinf{\und a}{\und q}{z}=\ftinf{\und a}{\und q
}{z}e^{(\infty)}_{\jtinf{\und a}}(z)$
with $\jtinf{\und a} = \text{diag}(1/\und{a})$.\\

\textit{Generators of the local Galois group at $\infty$.} We have two generators :
$$\breve{P}(y_0)^{-1}\text{diag}(e^{2\pi i\und\alpha})\breve{P}(y_0)
\text{ and } \breve{P}(y_0)^{-1}\text{diag}(\und u)\breve{P}(y_0).$$\\

\textit{Connection matrix}. The connection matrix is the limit as
$\und b$ tends to $\und q$ of :
\begin{eqnarray*}
(e^{(\infty)}_{\jtinf{\und a}}(z))^{-1} Q e^{(0)}_{\jtz{\und
q}}(z)\end{eqnarray*} \normalsize with :\scriptsize
$$Q=\pmatricet {\frac{\pochamer{\und
a_{\check1},\und{b}_{\check1}/a_1}{q}{\infty}}{\pochamer{\und
b_{\check1},\und a_{\check1}/a_1}{q}{\infty}} \frac{\theta_q (a_1
z)}{\theta_q(z)}} {\frac{\pochamer{\frac{q}{b_2} \und
a_{\check1},\und
b_{\check2}/a_1}{q}{\infty}}{\pochamer{\frac{q}{b_2} \und
b_{\check2},\und a _{\check 1}/a_1}{q}{\infty}} \frac{\theta_q
(\frac{qa_1}{b_2} z)}{\theta_q(z)}} {\frac{\pochamer{\frac{q}{b_3}
\und a_{\check1},\und b _{\check
3}/a_1}{q}{\infty}}{\pochamer{\frac{q}{b_3}\und b_{\check3},\und a
_{\check 1}/a_1}{q}{\infty}} \frac{\theta_q (\frac{qa_1}{b_3}
z)}{\theta_q(z)}} {\frac{\pochamer{\und
a_{\check2},\und{b}_{\check1}/{a_2}}{q}{\infty}}{\pochamer{\und
b_{\check1},\und a_{\check2}/a_2}{q}{\infty}} \frac{\theta_q (a_2
z)}{\theta_q(z)}} {\frac{\pochamer{\frac{q}{b_2} \und
a_{\check2},\und b_{\check
2}/a_2}{q}{\infty}}{\pochamer{\frac{q}{b_2}\und b_{\check 2},\und
a _{\check 2}/a_2}{q}{\infty}} \frac{\theta_q (\frac{qa_2}{b_2}
z)}{\theta_q(z)}} {\frac{\pochamer{\frac{q}{b_3} \und
a_{\check2},\und b_{\check
3}/a_2}{q}{\infty}}{\pochamer{\frac{q}{b_3}\und b _{\check 3},\und
a _{\check 2}/a_2}{q}{\infty}} \frac{\theta_q (\frac{qa_2}{b_3}
z)}{\theta_q(z)}} {\frac{\pochamer{\und
a_{\check3},\und{b}_{\check1}/a_3}{q}{\infty}}{\pochamer{\und
b_{\check1},\und a_{\check3}/a_3}{q}{\infty}} \frac{\theta_q (a_3
z)}{\theta_q(z)}} {\frac{\pochamer{\frac{q}{b_2} \und
a_{\check3},\und
b_{\check2}/a_3}{q}{\infty}}{\pochamer{\frac{q}{b_2} \und
b_{\check2},\und a _{\check 3}/a_3}{q}{\infty}} \frac{\theta_q
(\frac{qa_3}{b_2} z)}{\theta_q(z)}} {\frac{\pochamer{\frac{q}{b_3}
\und a_{\check3},\und b_{\check3}/a_3}{q}{\infty}}
{\pochamer{\frac{q}{b_3}\und b_{\check3},\und a _{\check
3}/a_3}{q}{\infty}} \frac{\theta_q (\frac{qa_3}{b_3}
z)}{\theta_q(z)}}
\pmatricet{1}{\frac{1-q^2/(b_2b_3)}{(\frac{q}{b_2}-1)(\frac{q}{b_3}-1)}}{\frac{q^2/(b_2b_3)}{(\frac{q}{b_2}-1)(\frac{q}{b_3}-1)}}
{0}{\frac{\frac{q}{b_3}-1}{(\frac{q}{b_3}-\frac{q}{b_2})(\frac{q}{b_2}-1)}}{\frac{-\frac{q}{b_3}}{(\frac{q}{b_3}-\frac{q}{b_2})(\frac{q}{b_2}-1)}}
{0}{\frac{\frac{q}{b_2}-1}{(\frac{q}{b_2}-\frac{q}{b_3})(\frac{q}{b_3}-1)}}{\frac{-\frac{q}{b_2}}{(\frac{q}{b_2}-\frac{q}{b_3})(\frac{q}{b_3}-1)}}$$
\normalsize which has the following form : \scriptsize
\begin{eqnarray*}
(e^{(\infty)}_{\jtinf{\und a}}(z))^{-1} \pmatricet {p_{1,1}(\und
a) \frac{\theta_q (a_1 z)}{\theta_q(z)}} {A_1\frac{\theta_q (a_1
z)}{\theta_q(z)}+p_{1,2}(\und a)a_1z\frac{\theta'_q (a_1
z)}{\theta_q(z)}} {B_2 \frac{\theta_q (a_1 z)}{\theta_q(z)} +
C_1z\frac{\theta'_q (a_1 z)}{\theta_q(z)} +
\frac{a_1^2}{2}p_{1,3}(\und a)z^2\frac{\theta''_q (a_1
z)}{\theta_q(z)}} {p_{2,1}(\und a) \frac{\theta_q (a_2
z)}{\theta_q(z)}} {A_2\frac{\theta_q (a_2
z)}{\theta_q(z)}+p_{2,2}(\und a)a_2z\frac{\theta'_q (a_2
z)}{\theta_q(z)}} {B_3 \frac{\theta_q (a_2 z)}{\theta_q(z)} +
C_2z\frac{\theta'_q (a_2 z)}{\theta_q(z)} +
\frac{a_2^2}{2}p_{2,3}(\und a)z^2\frac{\theta''_q (a_2
z)}{\theta_q(z)}} {p_{3,1}(\und a) \frac{\theta_q (a_3
z)}{\theta_q(z)}} {A_3\frac{\theta_q (a_3
z)}{\theta_q(z)}+p_{3,2}(\und a)a_3z\frac{\theta'_q (a_3
z)}{\theta_q(z)}} {B_3 \frac{\theta_q (a_3 z)}{\theta_q(z)} +
C_3z\frac{\theta'_q (a_3 z)}{\theta_q(z)} +
\frac{a_3^2}{2}p_{3,3}(\und a)z^2\frac{\theta''_q (a_3
z)}{\theta_q(z)}} e^{(0)}_{\jtz{\und
q}}(z)\end{eqnarray*}\normalsize where :

\begin{eqnarray*} p_{i,j}=p_{i,j}(\und a)= \frac{\pochamer{\und
a_{\check i},\und{q}_{\check j}/a_i}{q}{\infty}}{\pochamer{\und q
_{\check j},\und a_{\check i}/a_i}{q}{\infty}}.
\end{eqnarray*}
\\

\textit{Twisted connection matrix.}
\begin{scriptsize} \begin{eqnarray*}
\breve{P}(z)&=&\text{diag}((1/z)^{-\und \alpha}) \cdot \\ && \ \ \ \cdot \pmatricet
{p_{1,1}(\und a) \frac{\theta_q (a_1 z)}{\theta_q(z)}}
{A_1\frac{\theta_q (a_1 z)}{\theta_q(z)}+p_{1,2}(\und
a)a_1z\frac{\theta'_q (a_1 z)}{\theta_q(z)}} {B_2 \frac{\theta_q
(a_1 z)}{\theta_q(z)} + C_1z\frac{\theta'_q (a_1 z)}{\theta_q(z)}
+ \frac{a_1^2}{2}p_{1,3}(\und a)z^2\frac{\theta''_q (a_1
z)}{\theta_q(z)}} {p_{2,1}(\und a) \frac{\theta_q (a_2
z)}{\theta_q(z)}} {A_2\frac{\theta_q (a_2
z)}{\theta_q(z)}+p_{2,2}(\und a)a_2z\frac{\theta'_q (a_2
z)}{\theta_q(z)}} {B_3 \frac{\theta_q (a_2 z)}{\theta_q(z)} +
C_2z\frac{\theta'_q (a_2 z)}{\theta_q(z)} +
\frac{a_2^2}{2}p_{2,3}(\und a)z^2\frac{\theta''_q (a_2
z)}{\theta_q(z)}} {p_{3,1}(\und a) \frac{\theta_q (a_3
z)}{\theta_q(z)}} {A_3\frac{\theta_q (a_3
z)}{\theta_q(z)}+p_{3,2}(\und a)a_3z\frac{\theta'_q (a_3
z)}{\theta_q(z)}} {B_3 \frac{\theta_q (a_3 z)}{\theta_q(z)} +
C_3z\frac{\theta'_q (a_3 z)}{\theta_q(z)} +
\frac{a_3^2}{2}p_{3,3}(\und a)z^2\frac{\theta''_q (a_3
z)}{\theta_q(z)}} \cdot \\ && \ \ \ \ \ \ \ \ \ \ \ \ \ \ \ \ \ \ \ \ \ \ \ \ \ \ \ \ \ \ \ \ \ \ \ \ \ \ \ \ \ \ \ \ \ \ \ \ \ \ \ \ \ \ \ \ \ \ \ \ \ \ \ \ \ \ \ \ \ \ \ \ \ \ \ \ \ \ \ \ \ \ \ \ \ \ \ \ \ \ \ \ \ \ \ \ \ \ \ \ \ \ \ \ \ \ \ \ \ \ \ \ \ \ \cdot
\pmatricet{1}{\ell_q(z)}{\frac{\ell_q(z)(\ell_q(z)-1)}{2}}{0}{1}{\ell_q(z)}{0}{0}{1}
\end{eqnarray*}\end{scriptsize}

$_{}$\vskip 10 pt

\begin{theo}\label{lie irred bis}
Suppose that the $q$-hypergeometric equation $\hyp{\und a}{\und
q}$ is Lie-irreducible and has $q$-real
parameters. Then $\gzder=\sl{3}$. Moreover, we have :
\begin{itemize}
  \item[$\bullet$] $\G=\gl{3}$ if $a_1a_2a_3 \not \in
  q^\rinteger$;
  \item[$\bullet$] $\G=\overline{\langle \sl{3}, e^{2\pi i (\beta_2+\beta_3)}, v_1v_2 \rangle}$ if $a_1a_2a_3 \in q^\rinteger$.
\end{itemize}
\end{theo}

\begin{proof}
Arguing as for the proof of Theorem \ref{lie irred prem} we get
that $\gzder$ is either conjugated to $\psl$, or equal to
$\sl{3}$.

Suppose that $\gzder$ is conjugated to $\psl$ : there exists $R
\in \gl{3}$ such that $\gzder=R^{-1} \psl R$. The Galois group
$\G$ being a subgroup of the normalizer of $\gzder$, we deduce
from Lemma \ref{norm psl} that $\G \subset R^{-1} \pgl R$.

We have $\sharp \{ e^{2\pi i \alpha_1}, e^{2\pi i \alpha_2},
e^{2\pi i \alpha_3}\}=3$. Considering the generator of the local
Galois group at $0$ given above, we see that
$M=R\pmatricet{1}{1}{0}{0}{1}{1}{0}{0}{1}R^{-1}$ belongs to
$\pgl$. Let $T'$ be an invertible upper-triangular matrix such
that
$T'\pmatricet{1}{2}{1}{0}{1}{1}{0}{0}{1}T'^{-1}=\pmatricet{1}{1}{0}{0}{1}{1}{0}{0}{1}$.
We have $M=RT'\pmatricet{1}{2}{1}{0}{1}{1}{0}{0}{1}(RT')^{-1}\in
\psl$. On the other hand, $M$ being a non-semi-simple element of
$\pgl$, we deduce from Lemma \ref{lem diag psl} that there exists
$S \in \psl$ such that
$M=S\pmatricet{1}{2}{1}{0}{1}{1}{0}{0}{1}S^{-1}$. Considering the
two expressions of $M$ given above and using Lemma \ref{lem2}, we
see that there exists $T''$ an invertible upper-triangular matrix
such that $RT'=ST''$. Therefore the matrix $T=T''T'^{-1}$ is an
invertible upper-triangular matrix and satisfies $R=ST$.

Proceeding as in the proof of Theorem \ref{lie irred prem}, we see
that, for all $z\in
\bon$, there exist $\Sigma_z \in \permc$ 
such that $\Sigma_z \breve{P}(z)R^{-1} \in \psl$.

We conclude that, for all $z\in\bon$, $\Sigma_z
\breve{P}(z)T^{-1} \in \psl$.

Using the algebraic relations verified by the entries of the
elements of $\psl$, and arguing as for the proof of Theorem
\ref{lie irred prem}, we get that for some $i\in\{1,2,3\}$, a
functional equation of the following form holds on $\complex^*$ :

$$\left(*\frac{\theta_q (a_i z)}{\theta_q(z)}+p_{i,2}(\und
a)a_iz\frac{\theta'_q (a_i z)}{\theta_q(z)}\right)^2=K
\frac{\theta_q (a_i z)}{\theta_q(z)}\left(* \frac{\theta_q (a_i
z)}{\theta_q(z)} + *\frac{\theta'_q (a_i z)}{\theta_q(z)} +
\frac{a_i^2}{2}p_{i,3}(\und a)z^2\frac{\theta''_q (a_i
z)}{\theta_q(z)}\right)$$ for some constant $K\in \complex^*$ and
where each $*$ denotes a holomorphic function over $\complex^*
\setminus q^\real$. Now we get a contradiction : the right side of
the above equality vanishes for $z\in \frac{1}{a_i} q^\rinteger$
but not the left hand side. Indeed, the vanishing of the right
hand side on $z\in \frac{1}{a_i} q^\rinteger$ is clear; the
non-vanishing of the left hand side on $z\in \frac{1}{a_i}
q^\rinteger$ is a consequence of the fact that $\theta_q$ vanishes
exactly to the order one on $q^\rinteger$ and that $p_{i,2}(\und
a)\neq 0$ (because the equation is irreducible). We get a
contradiction.

Finally : $\gzder=\sl{3}$.

It is now easy to complete the proof using formula (\ref{det
connec}).
\end{proof}


\subsection{$\hbox{\underline a =(a,a,a) \text{ and } \und{b}=\und{q}}$}

\textbf{\textit{Hypotheses.}} In this section we suppose that :
$\und a=(a,a,a)$ and $\und b=\und q$.\\

\noindent \textit{\textbf{Basic objects.}}

\textit{Local solution at $0$.} We are in the non-resonant
logarithmic case at $0$ : the situation is the same than in the
previous section. The matrix valued function $\ytz{a\und 1}{\und
q}{z}=\ftz{a\und 1}{\und q}{z}e^{(0)}_{\jtz{\und q}}(z)$ is a
fundamental system of solutions of the
$q$-hypergeometric equation with $\und a =a \und 1$ and $\und b=\und q$.\\
\\

\textit{Generators of the local Galois group at $0$.} We have the
following generator :
$$\pmatricet{1}{1}{0}{0}{1}{1}{0}{0}{1}.$$\\

\textit{Local solution at $\infty$.} We have :
$$\qhypermatricetcour{\und a}{\und
b}{\infty}=\pmatricet{1/a^2}{-1/a}{1}{1/a^3}{0}{0}{1/a^4}{1/a^3}{0}
\jtinf{a \und 1}
\pmatricet{1/a^2}{-1/a}{1}{1/a^3}{0}{0}{1/a^4}{1/a^3}{0}^{-1}$$
with $\jtinf{a \und 1} =
\pmatricet{1/a}{1}{0}{0}{1/a}{1}{0}{0}{1/a}$. We are in the
non-resonant logarithmic case at $\infty$. A fundamental system of
solutions at $\infty$ of (\ref{syst hypergeo}) as described in
section \ref{section the basic objects} is given by $\ytinf{a \und
1}{\und q}{z}=\ftinf{a \und 1}{\und q}{z}e^{(\infty)}_{\jtinf{\und
a}}(z)$ where $\ftinf{a \und 1}{\und q}{z}$ is the limit as $\und
a$ tends to $a\und 1$ of :
\begin{eqnarray*}&&\ftinf{\und a}{\und q}{z} V(q\und a)^{-1} \pmatricet{1/a^2}{-1/a}{1}{1/a^3}{0}{0}{1/a^4}{1/a^3}{0}\\&=&\ftinf{\und a}{\und q}{z}  \pmatricet{a_1^2}{a_1^3/a_2}{a_1^3/a_2}{0}{-a_1^2/a_2+a_1^3/a_2^2}{-a_1^2/a_3+a_1^3/a_3^2}{0}{1-2a_1/a_2+a_1^2/a_2^2}{1-2a_1/a_3+a_1^2/a_3^2}^{-1}.\end{eqnarray*}

\textit{Generators of the local Galois group at $\infty$.} We have
two generators :
$$e^{2\pi i  \alpha} I_3,\ u I_3
\text{ and } \breve{P}(y_0)^{-1}\pmatricet{1}{a}{0}{0}{1}{ a}{0}{0}{1}\breve{P}(y_0).$$\\

\textit{Connection matrix}. The connection matrix is the limit as
$\und a$ tends to $a \und 1$ of : \scriptsize
\begin{eqnarray*}
&&(e^{(\infty)}_{\jtinf{a\und 1}}(z))^{-1}
\pmatricet{a_1^2}{a_1^3/a_2}{a_1^3/a_2}{0}{-a_1^2/a_2+a_1^3/a_2^2}{-a_1^2/a_3+a_1^3/a_3^2}{0}{1-2a_1/a_2+a_1^2/a_2^2}{1-2a_1/a_3+a_1^2/a_3^2}\cdot
\\ && \ \ \ \ \ \ \ \ \ \ \ \ \cdot \pmatricet {p_{1,1}(\und a) \frac{\theta_q (a_1
z)}{\theta_q(z)}} {A_1\frac{\theta_q (a_1
z)}{\theta_q(z)}+p_{1,2}(\und a)a_1z\frac{\theta'_q (a_1
z)}{\theta_q(z)}} {B_2 \frac{\theta_q (a_1 z)}{\theta_q(z)} +
C_1z\frac{\theta'_q (a_1 z)}{\theta_q(z)} +
\frac{a_1^2}{2}p_{1,3}(\und a)z^2\frac{\theta''_q (a_1
z)}{\theta_q(z)}} {p_{2,1}(\und a) \frac{\theta_q (a_2
z)}{\theta_q(z)}} {A_2\frac{\theta_q (a_2
z)}{\theta_q(z)}+p_{2,2}(\und a)a_2z\frac{\theta'_q (a_2
z)}{\theta_q(z)}} {B_3 \frac{\theta_q (a_2 z)}{\theta_q(z)} +
C_2z\frac{\theta'_q (a_2 z)}{\theta_q(z)} +
\frac{a_2^2}{2}p_{2,3}(\und a)z^2\frac{\theta''_q (a_2
z)}{\theta_q(z)}} {p_{3,1}(\und a) \frac{\theta_q (a_3
z)}{\theta_q(z)}} {A_3\frac{\theta_q (a_3
z)}{\theta_q(z)}+p_{3,2}(\und a)a_3z\frac{\theta'_q (a_3
z)}{\theta_q(z)}} {B_3 \frac{\theta_q (a_3 z)}{\theta_q(z)} +
C_3z\frac{\theta'_q (a_3 z)}{\theta_q(z)} +
\frac{a_3^2}{2}p_{3,3}(\und a)z^2\frac{\theta''_q (a_3
z)}{\theta_q(z)}} e^{(0)}_{\jtz{\und q}}(z)\end{eqnarray*}
\normalsize which has the following form :

\begin{eqnarray*}
(e^{(\infty)}_{\jtinf{a \und 1}}(z))^{-1} \pmatricet {*}{*}{*}{\text{*}} {*} {*}
{\frac{\theta_q(a_1)^2}{\pochamer{q}{q}{\infty}^4}
\frac{\theta_q (a_1 z)}{\theta_q(z)}} {*\frac{\theta_q (a_1
z)}{\theta_q(z)}+\frac{\theta_q(a_1)^2}{\pochamer{q}{q}{\infty}^4}a_1z\frac{\theta'_q
(a_1 z)}{\theta_q(z)}}{*}e^{(0)}_{\jtz{\und
q}}(z)\end{eqnarray*}
where each $*$ denotes an holomorphic function over $\complex^* \setminus q^\real$.\\

\textit{Twisted connection matrix.}
\begin{small} \begin{eqnarray*}
\breve{P}(z)&=&(1/z)^{-\alpha}\pmatricet{1}{a\ell_q(z)}{a^{2}\frac{\ell_q(z)(\ell_q(z)-1)}{2}}{0}{1}{a\ell_q(z)}{0}{0}{1}^{-1} \cdot \\ && \ \ \ \ \ \ \ \ \ \ \ \ \ \ \ \ \ \ \ \ \ \ \ \ \ \cdot \pmatricet {*}{*}{*}{\text *}{*}{*}
{\frac{\theta_q(a_1)^2}{\pochamer{q}{q}{\infty}^4}
\frac{\theta_q (a_1 z)}{\theta_q(z)}} {*\frac{\theta_q (a_1
z)}{\theta_q(z)}+\frac{\theta_q(a_1)^2}{\pochamer{q}{q}{\infty}^4}a_1z\frac{\theta'_q
(a_1 z)}{\theta_q(z)}}{*}
\pmatricet{1}{\ell_q(z)}{\frac{\ell_q(z)(\ell_q(z)-1)}{2}}{0}{1}{\ell_q(z)}{0}{0}{1}
\end{eqnarray*}\end{small}
where each $*$ denotes an holomorphic function  over $\complex^*
\setminus q^\real$. $_{}$\vskip 10 pt

\begin{theo}
Let us suppose that the $q$-hypergeometric equation $\hyp{a \und
1}{\und q}$ is Lie-irreducible and has
$q$-real parameters. Then $\gzder=\sl{3}$. More precisely :
\begin{itemize}
  \item[$\bullet$] $\G=\gl{3}$ if $a^3 \not \in
  q^\rinteger$;
  \item[$\bullet$] $\G=\overline{\langle \sl{3}, e^{2\pi i (\beta_2+\beta_3)}, v_1v_2 \rangle}$ if $a^3 \in q^\rinteger$.
\end{itemize}
\end{theo}

\begin{proof}
Arguing as for the proof of Theorem \ref{lie irred prem} we get
that $\gzder$ is either conjugated to $\psl$ or equal to $\sl{3}$.

Suppose that $\gzder$ is conjugated to $\psl$ : there exists $R
\in \gl{3}$ such that $\gzder=R^{-1} \psl R$. The Galois group
$\G$ being included in the normalizer of $\gzder$, we deduce from
Lemma \ref{norm psl} that $\G \subset R^{-1} \pgl R$.

Arguing as for Theorem \ref{lie irred bis}, we get an invertible
upper-triangular matrix $T$ such that $RT^ {-1} \in \psl$.

For similar reasons,  for all $z\in \bon$ there exists $T_z$ an
invertible upper-triangular matrix such that
$T_z\breve{P}(z)R^{-1} \in \psl$.

We deduce that $T_z\breve{P}(z)T^{-1} \in \psl$.

Combining this fact with the explicit form of the connection
matrix given above, we get a functional equation of the form :
$$\frac{\theta_q(a_1)^2}{\pochamer{q}{q}{\infty}^4} \theta_q (a_1 z) * = (*\theta_q (a_1 z)+\frac{\theta_q(a_1)^2}{\pochamer{q}{q}{\infty}^4}a_1z\theta'_q (a_1 z))^2$$
where each $*$ denotes some holomorphic function over $\complex^*
\setminus q^\real$. Now we get a contradiction : the left hand
side of the above functional equation vanishes for $z=1/a_1$ but
this is not the case of the right hand side since
$\theta_q(a_1)\neq 0$ (the equation is irreducible).

Finally : $\gzder=\sl{3}$. The theorem follows easily from
(\ref{det connec}).

\end{proof}


\subsection{The remaining non-resonant cases}

We have treated the following non resonant cases :

\begin{itemize}
  \item[(i)] $\forall i\neq j,\ \ a_i/a_j \not \in q^\rinteger$ and $b_2/b_3,b_2,b_3 \not \in
  q^\rinteger$;
  \item[(ii)] $\forall i\neq j,\ \ a_i/a_j \not \in q^\rinteger$ and $b_2=b_3 \not \in
  q^\rinteger$;
  \item[(iii)] $\forall i\neq j,\ \ a_i/a_j \not \in q^\rinteger$ and
  $b_2=b_3=q$;
  \item[(iv)] $\und a=(a,a,a)$ and $\und b=\und q$.
\end{itemize}

We have shown that in each cases (in the $q$-real case) the Galois
group $\G$ has $\gzder=\sl{3}$. We claim that the same conclusion
holds for every non-resonant
Lie-irreducible generalized $q$-hypergeometric equation of order
three. Indeed, the reader will easily adapt the proofs of the
cases (i), (ii), (iii) or (iv) to the remaining non-treated
non-resonant ($q$-real) cases. For instance :
\begin{itemize}
  \item[$\bullet$] if $a_1=a_2=a_3$ and $b_2/b_3,b_2,b_3 \not \in
  q^\rinteger$, we proceed as for the case (iii);
  \item[$\bullet$] if $a_1=a_2\not \in q^\rinteger a_3$, we proceed as for the case (ii), etc.
\end{itemize}

\nocite{roquesqgalois}

\subsection{Proof of the main theorem}

\begin{theo}
Let $\G$ be the Galois group of a Lie-irreducible generalized
$q$-hypergeometric equation $\hyp{\und a}{\und b}$ of order three
with $q$-real parameters. Then $\gzder=\sl{3}$. More precisely :
\begin{itemize}
  \item[$\bullet$] $\G=\gl{3}$ if $\frac{a_1a_2a_3}{b_2b_3} \not \in
  q^\rinteger$;
  \item[$\bullet$] $\G=\overline{\langle \sl{3}, e^{2\pi i (\beta_2+\beta_3)}, v_1v_2 \rangle}$ if $\frac{a_1a_2a_3}{b_2b_3} \in q^\rinteger$.
\end{itemize}
\end{theo}

\begin{proof}
We have already shown that any non-resonant Lie-irreducible generalized $q$-hypergeometric equation of
order three with $q$-real parameters has a Galois group $\G$ such
that $\gzder=\sl{3}$. The proof follows easily from this and from
Proposition \ref{qhyp irr}, which allows us to reduce the problem
to the non-resonant cases.
\end{proof}

As a concluding remark, it would be interesting to understand what
happens to the difference Galois groups of the Lie-irreducible
generalized $q$-hypergeometric equations under consideration in
the present paper as $q$ tends to $1$, or, more generally, as
$|q|$ tends to $1$.

\bibliography{biblio}
\bibliographystyle{plain}

\noindent \textsc{$_{}$\ \ \ \ \ Julien Roques\\
$_{}$\ \ \ \ \ D\'epartement de Math\'ematiques et Applications\\
$_{}$\ \ \ \ \ \'Ecole Normale Sup\'erieure\\
$_{}$\ \ \ \ \ 45, rue d'Ulm\\
$_{}$\ \ \ \ \ 75230 Paris Cedex 05 - France}\\
$_{}$\ \ \ \ \ \ E-mail : \textsf{julien.roques@ens.fr}

\end{document}